\documentclass[11pt]{article}

\usepackage{amsfonts}
\usepackage{amssymb}
\usepackage{amsmath}
\usepackage{amsthm}
\usepackage{latexsym}
\usepackage{verbatim}
\usepackage{epsfig}
\usepackage{amsbsy}
\usepackage{amscd}
\usepackage{bm}

\usepackage{graphicx}

\usepackage{ifthen} 
\usepackage{xcolor}
\usepackage[colorlinks=true, linkcolor=teal, anchorcolor=teal,
            citecolor=blue,  filecolor=blue, menucolor=blue, pagecolor=teal,											urlcolor=blue, plainpages=false, pdfpagelabels]{hyperref}

\newcommand{\mymail}[1]{\href{mailto:#1}{\texttt{#1}}}

\usepackage{fancyhdr}
\usepackage{currfile}
\fancypagestyle{firststyle}{

\cfoot{{\footnotesize \texttt{\currfilename}, version: \texttt{\today} }} 
}

\newcommand{\setauthA}[1]{\def\authA{#1}}

\def\printA{\begin{tabular}{l} \authA \end{tabular}}

\newcommand{\makemytitle}[1]{\begin{center}{\textsf{\LARGE #1}}
  \end{center}
}

\usepackage[sf,bf]{titlesec}

\usepackage{algorithmic}
\usepackage{algorithm}

\usepackage[nonamebreak,square,numbers,sort]{natbib}

\usepackage[letterpaper, textheight = 676pt]{geometry}
\providecommand{\mc}[1]{\mathcal#1}
\providecommand{\mc}[1]{\mathcal#1}

\newcommand{\R}{{\mathbb R}}

\DeclareMathOperator{\E}{\mathbf{E}}
\DeclareMathOperator{\p}{\mathbf{P}}

\DeclareMathOperator{\tr}{tr}

\providecommand{\T}{\top} 

\providecommand{\wt}[1]{\widetilde{#1}}
\providecommand{\wh}[1]{\widehat{#1}}

\providecommand{\nnorm}[1]{ \lVert#1 \rVert}
\newcommand{\scp}[2]{\left\langle#1, #2\right\rangle}
\newcommand{\nscp}[2]{\langle#1, #2\rangle}

\newcommand{\blanco}[1]{  }

\newcommand{\deriv}[3]{%
\ifthenelse{#1 = 1}{\frac{d\,#2}{d\,#3}}{\frac{d^{{#1}} #2}{d{#3}^{{#1}}}}
}

\newcommand{\partials}[3]{%
\ifthenelse{#1 = 1}{\frac{\partial\,#2}{\partial\,#3}}{\frac{\partial^{#1}
    #2}{\partial#3^{#1}}}
} 

\def\su{\sum_{i=1}^n}

\def \coloneq{\mathrel{\mathop:}=}

\def \eps{\varepsilon}

\newtheorem{theo}{Theorem}
\newtheorem{propo}{Theorem}
\newtheorem{definitio}{Theorem}

\newtheorem{lemmaA}{Theorem}[section]

\newtheorem{defn}[definitio]{Definition}

\newtheorem{prop}[propo]{Proposition}

\newtheorem{lemmaApp}[lemmaA]{Lemma}

\newenvironment{bew}{\begin{proof}[Proof]}{\end{proof}}

\usepackage{subfigure}

\newboolean{isjournal}
\setboolean{isjournal}{true}

\newboolean{nocolors}
\setboolean{nocolors}{false}

\newboolean{usemtpro}
\setboolean{usemtpro}{false}

\ifthenelse{\boolean{nocolors}}{\hypersetup{colorlinks=false}}{}
\ifthenelse{\boolean{usemtpro}}{}{\newcommand{\wwhat}[1]{\widehat{#1}}}

\makeatletter
\newcommand\footnoteref[1]{\protected@xdef\@thefnmark{\ref{#1}}\@footnotemark}
\makeatother

\setauthA{{\bfseries Martin Slawski} \\ Department of Statistics \\ George Mason University \\ Fairfax, VA 22030, USA \\ \mymail{mslawsk3@gmu.edu}}

\def\dwn{d \wedge n}

\begin{document}
\thispagestyle{firststyle}

\makemytitle{{\bfseries {On Principal Components Regression, Random Projections, and Column Subsampling
    }}}
\vskip 3.5ex
{\large\begin{center}
\printA
\end{center}}

\vskip 3.5ex

\begin{abstract} \noindent Principal Components Regression (PCR) is a traditional tool for dimension reduction in linear regression that has been both criticized and defended. One concern about PCR is that obtaining the leading principal components tends to be computationally demanding for large data sets.  While random projections do not possess the optimality properties of the leading principal subspace, they are computationally appealing and hence have become increasingly popular in recent years. In this paper, we present an analysis showing that for random projections satisfying a Johnson-Lindenstrauss embedding property, the prediction error in subsequent regression is close to that of PCR, at the expense of requiring a slightly large number of random projections than principal components. Column sub-sampling constitutes an even cheaper way of randomized dimension reduction outside the class of Johnson-Lindenstrauss transforms. We provide numerical results based on synthetic and real data as well as basic theory revealing differences and commonalities in terms of statistical performance.   
\end{abstract}

\section{Introduction}

Principal Components Regression (PCR), first introduced in \cite{Kendall1957, Hotelling1933}, is perhaps the most basic approach to dimension reduction in linear regression. In PCR, the design matrix $X \in \R^{n \times d}$ containing the original predictor variables is replaced by $X R \in \R^{n \times r}$, $r < \dwn$, where $R \in \R^{r \times d}$ reduces $X$ to its top $r$ principal components. From a statistical point of view, PCR can be motivated as a way of dealing with multi-collinearity and reducing estimation variance at the expense of additional bias. From the point of view of computation, PCR potentially achieves a reduction from a large number of variables to a parsimonious model, which can be beneficial for both model fitting and prediction of future observations. The use of PCR is debated in the literature as it does not need to be case that principal components corresponding to small singular values do not significantly contribute in predicting the response variable \cite{Joliffe1982, Artemiou2009}. Herein, we mostly avoid touching upon this issue. Instead, the purpose of this paper is to establish a connection between PCR and the use of randomized methods of dimension reduction in linear regression in which the matrix $R$ above is sampled from a suitable distribution. The latter approach, typically referred to as the method of random projections, is motivated from large-scale datasets in which both the number of samples $n$ and the dimension $d$ are large; in this case, computation of principal components via the SVD can be demanding. Random projections only require a matrix multiplication which can be easily parallelized. Having its roots in the celebrated Johnson-Lindenstrauss-Lemma \cite{JohnsonLindenstrauss}, the idea has meanwhile a long history in computer science \cite{Vempala2005}, and has recently attracted considerable interest in statistics (see \cite{Cannings2017} and the references therein).
\vskip1ex
\noindent\emph{Contributions and Related Work.} It is critical to understand the statistical properties of the dimension reduction provided by random projections. Regarding linear regression, there are many more papers (e.g., \cite{RaskuttiMahoney2014, PilanciWainwright2015, Zhou2009, Wang2017, Homrighausen2017, Ahfock2017}) on the scenario in which $X$ is reduced to $RX$, i.e., $R$ is multiplied from left instead of from the right with $X$ being reduced to $XR$. The latter scenario was first analyzed in \cite{Maillard2009} under the term
``compressed least squares'' (CLS) which will also be employed herein. For a fixed design setting, refinements appear in \cite{Kaban2014}, and very recently in \cite{Thanei2017}. Together with a preliminary version \cite{Slawski2017} of the present paper, the paper \cite{Thanei2017} is the first to make a connection between PCR and CLS. However, as we show below, while improving over the main result in \cite{Kaban2014}, the upper bounds on the prediction error of CLS in \cite{Thanei2017} still leave a considerable gap to PCR. In the present paper, we try to close this gap. In brief, our main result states that CLS can roughly match the performance of PCR with respect to prediction at the expense of a moderate increase (at most by a logarithmic factor) in the reduced dimension. This property is shown to hold for a certain class of matrices comprising those that are typically considered in the literature on randomized dimensionality reduction and compressed sensing. We leave it as an open problem whether a similar result can be established for column subsampling, i.e.,~the columns of $R$ are chosen uniformly at random from the canonical basis vectors. Results are provided indicating that more stringent conditions are required in that case. Finally, we note that the very recent work \cite{Kasiviswanathan2017} presents an analysis of CLS when the goal is to recover the vector of regression coefficients under sparsity. Such an assumption is not made herein, and accordingly settings, goals and results are not comparable. 
\vskip2ex
\noindent \emph{Outline.} In Section $\S$2, we provide background and review the results of prior work in more detail. Our main result is contained in $\S$3, while $\S$4 discusses extensions and open questions. We conclude with a brief summary im $\S$5. The appendix contains all proofs.  

\vskip2ex
\noindent \emph{Notation.} For a positive integer $m$, we write $[m] \coloneq \{1,\ldots,m\}$. For a matrix $M$, we write $P_M$ for the projection operator onto the column space of $M$. Its Frobenius norm is denoted by $\nnorm{M}_F = \sqrt{\tr(M^{\T} M)}$, where ``$\tr$'' is the trace of a diagonal matrix. The Gaussian distribution with zero mean and variance $s^2$ is denoted by $N(0, s^2)$, and $\text{Unif}(S)$ denotes the uniform probability distribution on a set $S$. For $a, b \in \R$, we write $a \wedge b = \min\{a,b\}$ and
$a \vee b = \max\{a,b\}$. Positive constants are denoted by $C$, $C_1$, $c$, $c_1$ etc. We make use of the usual Big-O notation in terms of $O$, $o$, $\Omega$ and $\Theta$.




\section{Background}

We start by providing some context for our main result. After fixing the 
setup, we derive bounds on the prediction error of PCR and put them into relation to existing results on CLS. This will point to a significant gap that motivates our analysis of CLS in $\S$\ref{sec:cls_improved}. 

\vspace*{-0.1in}
\subsection{General setting}
\vspace*{-0.05in}
We consider fixed design linear regression for data  $(y_i, x_i)$, with $y_i$ taking values
in $\R$ and $x_i$ taking values in $\R^d$, $i \in [n]$. The predictors $x_i$ are considered as fixed, and
\begin{equation*}
y_i = f_i + \xi_i,
\end{equation*}
with $f_i = \E[y_i]$ and $\xi_i$ following a distribution with mean zero and variance $\sigma^2$, $i \in [n]$. Moreover, the $\{ \xi_i \}_{i = 1}^n$ are assumed to be uncorrelated. More concisely, we write $y = f + \xi$, where $y = (y_i)_{i=1}^n$, etc. We denote by $X \in \R^{n \times d}$ the design matrix whose rows are given by $x_i^{\T}$, $i \in [n]$. The optimal linear predictor $Xw^*$ of $y$ given $X$ with respect to squared loss is defined by the optimization problem
\begin{equation*}
\min_{w \in \R^d} \E[\nnorm{y - X w}_2^2/n],
\end{equation*}
where the expectation is with respect to the noise $\xi$. Any minimizer $w^*$ of the above problem satisfies $X w^* = P_X f$ with $P_X$ defined as in the paragraph on notation above;  
if there are multiple such $w^*$ we choose the one with minimum $\ell_2$-norm. Accordingly, we define the excess risk of an estimator $\wh{\theta} = \wh{\theta}(X,y)$ of $w^*$ by
\begin{equation*}
\mc{E}(\wh{\theta}) =  \E[\nnorm{Xw^* - X\wh{\theta}}_2^2/ n],
\end{equation*}
If the linear model holds exactly (i.e.,~$P_X f = f$), $\mc{E}(\wh{\theta})$ equals the in-sample mean squared prediction error that measures how well the 
$\{ x_i^{\T} \wh{\theta} \}_{i = 1}^n$  predict the ``denoised'' observations $\{ x_i^{\T} w^* \}_{i = 1}^n$ on average. An ordinary least squares (OLS) estimator $\wh{w}$ satisfies $X \wh{w} = P_X y$. Its excess risk is given by 
\begin{equation}\label{eq:excess_ols}
\mc{E}(\wh{w}) = \sigma^2 \text{rank}(X) / n. 
\end{equation}
To keep matters simple, we assume that $X$ has full rank $d \wedge n$ unless otherwise stated. In this paper, we are interested in a high-dimensional setup in which $\text{rank}(X)$ is of the same order of magnitude as $n$. In this situation, OLS does not yield satisfactory statistical performance. Moreover, if both $n$ and $d$ are large, obtaining $\wh{w}$ or making predictions based on $\wh{w}$ becomes computationally costly.

In light of these issues, it makes sense to consider alternatives that aim at leveraging some sort of low-dimensional structure. Scenarios in which $w^*$ exhibits one of various forms of sparsity are dominating in the literature, see the monographs of \cite{Buehlmann2011, HastieTibWain2015} for an overview. In the present paper, we follow another direction in which the predictors $\{ x_i \}_{i = 1}^n$ are linearly mapped into a lower-dimensional space, and linear least squares regression is then performed based on the subspace obtained in this way. Put differently, one considers a new design matrix $X_R = X R$ with $R$ being a $d$-by-$k$ matrix, $k \ll d$. On the statistical side, one potentially achieves a substantial reduction of variance at the expense of an increased bias as made precise below. The excess risk of the approach is given by   
\begin{equation}\label{eq:excess_R}
  \mc{E}(R) = \E \left[\nnorm{X w^* - X_R \wh{w}_R}_2^2 / n \right], 
\end{equation}
where $\wh{w}_R$ is a least squares solution based on the reduced design matrix $X_R$,
i.e., $\wh{w}_R$ satisfies $X_R \wh{w}_R = P_{X_R} y$, and the expectation is with respect to $\xi$ (in later sections, $R$
will be random, and we will then also take the expectation with respect to $R$).
Straightforward calculations show that $\mc{E}(R)$ can be decomposed into a bias and a variance term:
\begin{equation}\label{eq:excess_cls}
\mc{E}(R) = \underbrace{\nnorm{(I - P_{X_R}) Xw^*}_2^2/n}_{\text{Bias}} \; + \underbrace{\sigma^2 \text{rank}(X_R)/n}_{\text{Variance}}.  
\end{equation}
We commonly have $\text{rank}(X_R) = k$. The choice of $k$ determines the bias-variance trade-off. If $k$ can be chosen much smaller than $\dwn$ while at the same time the magnitude of the bias can be controlled, an improvement over the excess risk of OLS in \eqref{eq:excess_ols} is obtained.

The approach can as well be motivated from the computational side given that we only need to solve a least squares problem of dimension $k$ instead of $\dwn$. In addition, having a smaller number of predictors yields savings in storage and when making predictions.





\subsection{Excess risk of PCR}\label{sec:pcr}
The traditional choice of constructing $X_R$ is in terms of the leading principal components of $X$. Subsequent use of this new
set of predictors in regression is known as principal components regression (PCR). Let $X = U \Sigma V^{\T}$ be the singular
value decomposition (SVD) of $X$, where $U \in \R^{n \times \dwn}$, $U^{\T} U = I$, is the matrix of left singular vectors, $\Sigma \in \R^{\dwn \times \dwn}$ is the diagonal matrix whose diagonal contains the decreasingly ordered sequence of singular values $\sigma_1 \geq \ldots \geq \sigma_{\dwn}$, and
$V \in \R^{d \times \dwn}$, $V^{\T} V = I$, is the matrix of right singular vectors. For $r \in \{1,\ldots,\dwn\}$, consider
\begin{equation}\label{eq:partitioning}
\hspace*{-1ex}  U = [U_r \,\; U_{r+}], \qquad\, \Sigma = \begin{bmatrix}
    \Sigma_r &  0 \\
    0        & \Sigma_{r+}
  \end{bmatrix}, \qquad\,  V = [V_r \;\, V_{r+}],
\end{equation}
where $U_r$ and $V_r \in \R^{d \times r}$ contain the top $r$ left respectively right singular vectors, and
$\Sigma_r$ contains the corresponding singular values. The remaining singular vectors respectively singular values are contained in $U_{r+}$, $V_{r+}$ and $\Sigma_{r+}$. The top $r$ principal components are extracted from $X$ by setting  $R = V_r$:
\begin{align*}
  X_R = X V_r = (U_r \Sigma_r V_r^{\T} + U_{r+} \Sigma_{r+} V_{r+}^{\T}) V_r 
              = U_r \Sigma_r.      
\end{align*}
The corresponding projection is given by $P_{X_R} = U_r U_r^{\T}$ and the bias term in \eqref{eq:excess_cls} results as
\begin{align}\label{eq:P_pcr}
  (I - P_{X_R}) Xw^* = (I - U_r U_r^{\T}) Xw^* 
                    = U_{r+} \Sigma_{r+} V_{r+}^{\T} w^*.     
\end{align}
Let us define $\alpha^* \in \R^{\dwn}$ by
\begin{equation}\label{eq:alphastar}
  \alpha^* = V^{\T} w^* = \left[ \begin{array}{c}
                                   V_r^{\T} w^* \\
                                   V_{r+}^{\T} w^* 
                                 \end{array} \right] =  \left[ \begin{array}{c}
                                   \alpha_r^* \\
                                   \alpha_{r+}^{*} 
                                 \end{array} \right]. 
                               \end{equation}
Combining \eqref{eq:excess_R}, \eqref{eq:P_pcr} and \eqref{eq:alphastar}, the excess risk of PCR can then be expressed as follows.
\begin{align}
 \mc{E}(V_r) &= \nnorm{U_{r+} \Sigma_{r+} V_{r+}^{\T} w^*}_2^2/n + \sigma^2 r/n \notag\\
              &= \nnorm{\Sigma_{r+} \alpha_{r+}^*}_2^2/n +  \sigma^2 r/n \notag\\
              &=  \sum_{j=r+1}^{\dwn} \sigma_j^2 (\alpha_j^*)^2/n + \sigma^2 r / n \label{eq:excessrisk_pcr_exact}
\end{align}
We see from \eqref{eq:excessrisk_pcr_exact} that the excess risk of PCR behaves favorably if (i) the tail of the squared singular values at truncation level $r$ is small (i.e.,~$X$ can be well approximated by a matrix of rank $r$) and (ii) if there are no large coefficients in $\alpha^*$ outside its top $r$ entries corresponding to the leading singular vectors. Condition (ii) constitutes the main source of criticism of PCR: if nature is malicious, then $\alpha^*$ has most of its mass in $\alpha_{r+}^*$. Under generic random sampling, however, this is not a concern: if $V$ is sampled uniformly at random from its respective Stiefel manifold, then $\alpha^* / \nnorm{\alpha^*}_2$ is uniformly distributed on the unit sphere in $\R^{\dwn}$ so that the entries of $\alpha^*$ are roughly homogeneous in magnitude.   

In the sequel, we derive a series of bounds on the excess risk of PCR depending on the decay of the squared singular values
$\{ \sigma_j^2 \}_{j = 1}^{\dwn}$. For this purpose, we use the following simple upper bound on $\mc{E}(V_r)$:
\begin{equation}\label{eq:excessrisk_pcr_bound_dense}
  \mc{E}(V_r) \leq \nnorm{\alpha^*}_{\infty}^2 \frac{\nnorm{\Delta_r}_F^2}{n} + \sigma^2 \frac{r}{n}, \qquad \Delta_r \coloneq X - \mc{T}_r(X),
\end{equation}
where $\mc{T}_r(X) = U_r \Sigma_r V_r^{\T}$ equals the best rank $r$-approximation to $X$ with respect to the Frobenius norm. For what follows, we assume that $X$ is scaled such that $\nnorm{X}_F^2 = \sum_{j = 1}^{\dwn} \sigma_j^2 = n \cdot d$.  For $1 \leq s \leq \dwn$, we define
\begin{equation}
\gamma(s) = \sum_{j = 1}^s \sigma_j^2, \qquad \tau(s) = \frac{\gamma(\dwn) - \gamma(s)}{\gamma(\dwn)}.  
\end{equation}
The quantity $\nnorm{\Delta_r}_F^2 / n$ in \eqref{eq:excessrisk_pcr_bound_dense} can then be expressed as
\begin{align}\label{eq:gammatauDelta}
  \nnorm{\Delta_r}_F^2/n =\{\gamma(\dwn) - \gamma(r)\}/n = \tau(r) \cdot \gamma(\dwn) / n = \tau(r) \cdot d. 
\end{align}
After these preparations, we study the excess risk of PCR in three basic scenarios. 
\vskip1ex
\noindent \textbf{Scenario (F): perfectly flat spectrum}
\vskip1ex
\noindent
A flat spectrum means that $\sigma_j^2 = n \vee d$, $j \in [\dwn]$. We obtain that 
$\nnorm{\Delta_r}_F^2/n = d \cdot \tau(r) = \frac{d}{n} \vee 1 \, (\dwn - r)$. The choice of $r$ minimizing
the bound \eqref{eq:excessrisk_pcr_bound_dense} is given by $r^* = \dwn$ if $\nnorm{\alpha^*}_{\infty} > \sigma / \sqrt{d \vee n}$, and $r^* = 0$ otherwise. We note that when using the exact expression \eqref{eq:excessrisk_pcr_exact} for the excess risk, the optimal $r^*$ would result as the largest value of $r$ such that $\alpha_{r} > \sigma /  \sqrt{d \vee n}$. Eventually, this does not make much of a difference if the entries of $\alpha^*$ are of a comparable magnitude, and does not affect the conclusion that in general, we cannot hope for improvements over OLS when the spectrum is constant.   
\vskip1ex
\noindent \textbf{Scenario (P): polynomial decay}
\vskip1ex
\noindent
Suppose that $\sigma_j^2 = C \cdot j^{-q}$, $j \in [d]$, for $q \geq 2$ and a constant $C$ determined by the
relation $\sum_{j = 1}^{\dwn} \sigma_j^2 = n\cdot d$. Comparing series and integrals, we obtain $\gamma(\dwn) - \gamma(r) \leq C (q-1)^{-1} r^{-(q-1)}$. Moreover, $\gamma(\dwn) \geq C$, so that $\tau(r) \leq (q-1)^{-1} r^{-(q-1)}$. By 
\eqref{eq:excessrisk_pcr_bound_dense}
\begin{equation*}
\mc{E}(V_r) \leq (q-1)^{-1} r^{-(q-1)} \cdot d \cdot \nnorm{\alpha^*}_{\infty}^2 + \sigma^2 r/n.
\end{equation*}
Minimizing the right hand side w.r.t.~$r$, we obtain
\begin{align}\label{eq:excess_pcr_poly}
    r^* = \left\{ \nnorm{\alpha^*}_{\infty}^2 (n \cdot d) \big / \sigma^2   \right\}^{1/q}, \qquad
  \mc{E}(V_{r^*}) \leq 2  \left( d \nnorm{\alpha^*}_{\infty}^2 \right)^{1/q} \left( \sigma^2 / n \right)^{(q-1)/q}.
\end{align}
To get some insight into \eqref{eq:excess_pcr_poly}, fix $q = 2$ and consider the case of generic random sampling of $V$ as discussed
above so that $\alpha^* / \nnorm{\alpha^*}_2$ follows a uniform distribution on the unit sphere in $\R^{\dwn}$. In this situation,
$\nnorm{\alpha^*}_{\infty} / \nnorm{\alpha^*}_2$ scales as $O(\sqrt{\log(d) / d})$ as $d$ gets large. Assuming further
that $\nnorm{\alpha^*}_2 = O(1)$ yields $r^* = O(\{\log(d) n\}^{1/2})$ and $\mc{E}(V_{r^*}) = O(\sqrt{\log(d)/n})$. We have hence identified a regime in which PCR achieves better statistical and computational performance than OLS if $d = \Omega(n)$. Clearly, the improvements get amplified as $q$ increases.
\vskip1ex
\noindent \textbf{Scenario (E): exponential decay}
\vskip1ex
\noindent
Suppose that $\sigma_j^2 = C_0 \theta^j$ for $\theta \in (0,1)$. Then, $\tau(r) \leq \frac{\theta^r}{1-\theta}  = C_1 \exp(-c r)$, say. The optimal choice of $r^*$ and the corresponding bound on $\mc{E}(V_{r^*})$ result as
\begin{align}\label{eq:excess_pcr_exp}
  \begin{split}
  r^* = \frac{1}{c} \log\left(C_2 \nnorm{\alpha^*}_{\infty}^2 \, n  \,d  \big / \sigma^2 \right), \qquad
  \mc{E}(V_{r^*}) \leq \frac{2}{c} \left\{ \log\left(C_2 \nnorm{\alpha^*}_{\infty}^2  \, n  \,d  \big / \sigma^2 \right)  \, \vee 1  \right \}\,\sigma^2 / n. 
  \end{split}
\end{align}
\vskip2ex
\noindent Cases \textbf{(P)} and \textbf{(E)} show that PCR may improve significantly over OLS in terms of achievable dimension reduction and excess risk depending on the decay of the spectrum of $X$.

\subsection{Existing bounds for CLS}\label{sec:cls_existing}

We now consider the case of dimension reduction via a \emph{random} matrix $R$. We refer to the columns of $R$ as ``random projections'' as $R$ maps the predictors $\{ x_i \}_{i = 1}^n$ to a random linear subspace, typically of dimension $k$. Regarding the distribution of $R$, sampling its entries i.i.d.~from a Gaussian distribution with expectation zero and variance $1/k$ constitutes the basic case in the literature on randomized dimensionality reduction \cite{Vempala2005}. The column space of $R$ then follows the uniform distribution on the Grassmannian $\textsf{G}(d,k)$. Random Gaussian matrices of this form are the canonical example of Johnson-Lindenstrauss transforms \cite{JohnsonLindenstrauss} (henceforth JLTs for short), cf.~Definition \ref{defn:JLT_1} below. This class of matrices extends to i.i.d.~sub-Gaussian 
matrices \cite{Achlioptas2003, Matousek2008}, the fast JLT of \cite{Ailon2006}, and certain row-subsampled orthonormal matrices
\cite{AilonLiberty2011, KrahmerWard2011, Tropp2011} as they are also used in compressed sensing \cite{CandesWakin2008}.



Maillard \& Munos \cite{Maillard2009} were the first to study the use of JLTs for randomized
dimension reduction in least squares regression with random design under the name ``compressed least squares'' (CLS). They show a bound on a corresponding notion of excess risk of the order
\begin{equation}\label{eq:excess_MaillardMunos}
O\left(\sigma \nnorm{w^*}_2 \left(\E[\nnorm{x_1}_2^2] \right)^{1/2} \sqrt{\log(n) / n} \right)
\end{equation}
for $k = \Theta(\sqrt{n \log n} \nnorm{w^*}_2 \E[\nnorm{x_1}^2]^{1/2} / \sigma)$ random projections. For fixed design,
Kaban \cite{Kaban2014} (specializing Theorem 1 therein to the case where $R$ has i.i.d.~$N(0,1/k)$ entries) shows that
\begin{align}
  \E[\mc{E}(R)] &\leq \frac{(w^*)^{\T}(\tr(\Gamma) I + \Gamma) w^*}{k} + \sigma^2 \frac{k}{n},  \qquad \Gamma \coloneq X^{\T} X / n,\notag \\
                &\leq c  \frac{\tr(\Gamma) \nnorm{w^*}_2^2}{k}  + \sigma^2 \frac{k}{n}  = c \cdot \nnorm{w^*}_2^2 \left(\frac{1}{n} \sum_{j = 1}^{\dwn} \sigma_j^2 \right) +  \sigma^2 \frac{k}{n}\; \;\,\text{for some $c \in [1,2]$},  \label{eq:excess_Kaban_raw}
\end{align}
where $\mc{E}(R)$ is defined in \eqref{eq:excess_R} and the expectation in \eqref{eq:excess_Kaban_raw} is with respect to $R$.


Optimizing the bound \eqref{eq:excess_Kaban_raw} with respect to $k$, we obtain
\begin{equation}\label{eq:excess_Kaban_simple}
k^* = c' \sqrt{n \, \tr(\Gamma)}/ \sigma, \qquad \E[\mc{E}(R)] \leq 2 \sigma c'  \nnorm{w^*}_2 \sqrt{\tr(\Gamma)/n}. 
\end{equation}
with $c' = \sqrt{c}$. Comparing the bound of Maillard and Munos \eqref{eq:excess_MaillardMunos} with \eqref{eq:excess_Kaban_simple}, we essentially observe
an agreement apart from a $\sqrt{\log n}$ factor, noting that for random design with $\{x_i\}_{i = 1}^n \sim x$ such that $\E[xx^{\T}] = \Gamma_0$,
we have $\sqrt{\E[\nnorm{x_1}^2]} = \sqrt{\tr(\Gamma_0)}$.

In \cite{Maillard2009, ShahMeinshausen2013}, the excess risk bound in \eqref{eq:excess_Kaban_simple} is interpreted as  
being of the order $O(1/\sqrt{n})$. This would mean that the performance of CLS is comparable to that of PCR in scenario $(\textbf{P})$ with exponent
$q = 2$ above, independent of the $\{\sigma_j^2 \}_{j = 1}^{\dwn}$. However, such interpretation is not valid in general. In a fixed design setting,
it is common to assume that the columns $\{X_j \}_{j = 1}^d$ of $X$ are scaled such that $\nnorm{X_j}_2^2 = n$, $j \in [d]$, whereas for standard random designs, e.g.~$X$ with i.i.d.~rows from a zero-mean Gaussian distribution with unit variances, this scaling holds in expectation. In this case, $\tr(\Gamma)$ respectively $\E[\nnorm{x_1}_2^2]$ evaluate as $d$ which makes the bounds \eqref{eq:excess_MaillardMunos} and \eqref{eq:excess_Kaban_simple} of rather limited use. For \eqref{eq:excess_Kaban_simple}, we obtain
\begin{equation*}
k^* = c' \sqrt{n \, d}/ \sigma, \qquad \E[\mc{E}(R)] \leq 2 \sigma c'  \nnorm{w^*}_2 \sqrt{d/n}. 
\end{equation*}
This means that $k^*$ is of the same order (or may even exceed) $\dwn$, while the bound on the excess risk is generally inferior to that of OLS \eqref{eq:excess_ols}.

It turns out that the outcome \eqref{eq:excess_Kaban_simple} is the consequence of a crude bound on the bias of CLS. From \eqref{eq:excess_Kaban_raw}, we find that
the correct variance term $\sigma^2 k/n$ is present. On the other hand, a simple argument shows that the bias term $\tr(\Gamma) \nnorm{w^*}_2^2 / k$ is improvable: if
$X$ has rank $r$, $1 \leq r \leq \dwn$, CLS with $R$ as a matrix with $N(0,1)$-entries yields $P_X = P_{X_R}$ and in turn zero bias in \eqref{eq:excess_R} with probability one as long as $k \geq r$. By contrast, according to \eqref{eq:excess_Kaban_raw} basically $k/d \rightarrow \infty$ is required for the bias to vanish.

In \cite{Kaban2014}, the expected bias (where the expectation is w.r.t.~$R$) is bounded as
\begin{align}
  \E[\nnorm{(I - P_{X_R}) Xw^*}_2^2/n] &= \E \left[\min_{v \in \R^k} \nnorm{Xw^* - X R v}_2^2 / n\right] \leq \E [\nnorm{Xw^* - X R R^{\T} w^*}_2^2 / n] \label{eq:kaban_bias1}
\end{align}
Evaluation of \eqref{eq:kaban_bias1} then inevitably leads to  the term $\tr(\Gamma) \nnorm{w^*}_2^2 / k$ as it appears in \eqref{eq:excess_Kaban_raw}.

In a recent paper by Thanei et al. \cite{Thanei2017}, the following bound is used instead of \eqref{eq:kaban_bias1}:
\begin{equation}\label{eq:thanei_bias}
 \E \left[\min_{v \in \R^k} \nnorm{Xw^* - X R v}_2^2 / n\right] \leq \min_{v \in \R^k} \E[\nnorm{Xw^* - X R v}_2^2 / n] 
\end{equation}
The authors evaluate the above expectation for $R$ with i.i.d.~$N(0,1/k)$ entries, and then minimize with respect to
$v$. This yields the bound (cf.~Theorem 2 in \cite{Thanei2017})
\begin{align}
  &\E[\mc{E}(R)] \leq \frac{1}{n} \sum_{j = 1}^{\dwn} (\alpha_j^*)^2 \sigma_j^2 \omega_j + \sigma^2 \frac{k}{n}, \label{eq:bound_thanei}\\
  &\qquad \quad \omega_j = \frac{(1 + 1/k) \sigma_j^4 + (1 + 2/k) \sigma_j^2 \tr(\Gamma) + \tr(\Gamma)^2 / k}{(k + 2 + 1/k) \sigma_j^4 + 2(1 + 1/k) \sigma_j^2 \tr(\Gamma) + \tr(\Gamma)^2 / k}, \;\; j \in [\dwn]. \label{eq:weights_thanei}   
\end{align}
In order to assess \eqref{eq:bound_thanei} for improvements over earlier bounds, one needs to gain more insights into the weights $\{ \omega_{j} \}$. In Appendix \ref{app:bound_thanei_lower}, we show that $\min_{1 \leq j \leq \dwn} \omega_j \geq 2/(2+k)$, which implies that 
\begin{equation}\label{eq:bound_thanei_lower}
\sum_{j = 1}^{\dwn} (\alpha_j^*)^2 \sigma_j^2 \omega_j + \sigma^2 \frac{k}{n} \geq
  \frac{2}{2 + k} \frac{1}{n} \sum_{j = 1}^{\dwn} (\alpha_j^*)^2 \sigma_j^2 + \sigma^2 \frac{k}{n} =  \frac{2}{2 + k}  (w^*)^{\T} \Gamma w^* + \sigma^2 \frac{k}{n}.
\end{equation}
Compared to \eqref{eq:excess_Kaban_raw}, this means that at best, the term $\tr(\Gamma) \nnorm{w^*}_2^2$ gets replaced by $(w^*)^{\T} \Gamma w^*$. However,
even the lower bound in \eqref{eq:bound_thanei_lower} does not yield satisfactory results if $X$ is (approximately) of low rank as the bias term
again scales as $O(1/k)$ independent of the spectrum. As a consequence, no matter how small the rank of $X$ is, according to \eqref{eq:bound_thanei_lower} we still obtain an upper bound on the bias of $O(1/\sqrt{n})$ and $k^* = \Omega(\sqrt{n})$ for the optimal number of random projections.

\section{Improved Analysis (of CLS)}\label{sec:cls_improved}
In this section we present and discuss the main result of the paper, a bound on the excess risk of CLS that is of a similar flavor
of that of PCR in \S\ref{sec:pcr}. In this manner, we establish a substantially stronger connection between CLS and PCR 
as could be made based on other analyses reviewed in \S\ref{sec:cls_existing}.

\subsection{Assumptions on the random projections}
Our analysis basically requires $R$ to be a Johnson-Lindenstrauss transform (JLT). The precise conditions are
given as follows.
\begin{defn}\label{defn:JLT_1} Fix an arbitrary set of points of $m$ points $\{v_1,\ldots, v_m \} \subset \R^d$. A random $d$-by-$k$ matrix $R$ is said to be an $(m, \eps, \delta)$-JLT for $\eps, \delta \in (0,1)$ if
   \begin{equation*}
     (1 - \eps)\nnorm{v_i}_2^2 \leq \nnorm{R^{\T} v_i}_2^2 \leq  (1 + \eps)\nnorm{v_i}_2^2, \;\; i \in [m],
\end{equation*}
holds with probability at least $1 - \delta$. 
\end{defn}
The next definition is akin to the restricted isometry property  in the theory of sparse estimation \cite{Baraniuk2006}  with the difference that approximate norm preservation is required only for a single subspace (as opposed to the union of subspaces of sparse vectors).
\begin{defn} Let $\mc{V} \subset \R^d$ be an arbitrary subspace of dimension $s$. A random $d$-by-$k$ matrix $R$ is said to be an $(s,\eps,\delta)$-restricted
 isometry for $\eps, \delta \in (0,1)$ if 
 \begin{equation*}
 (1 - \eps)\nnorm{v}_2 \leq \nnorm{R^{\T} v}_2 \leq  (1 + \eps)\nnorm{v}_2 \;\;\, \forall v \in \mc{V}
\end{equation*}
holds with probability at least $1 - \delta$. 
\end{defn}
\noindent \textbf{Remark.} The conditions in the above two definitions are related in the following way: it is shown in~\cite{Baraniuk2006} that if $R$ is a $\left(\left\{\frac{12}{\eps}\right\}^k, \eps/2, \delta \right)$-JLT, then $R$ is also an $(s,  \eps, \delta)$-restricted isometry. We here state two separate definitions for ease of reference.  
\subsection{Main result}
We are now in position to state our main result. 
\begin{theo}\label{theo:mainresult} For $r \in \{1,\ldots,\dwn\}$, let the following conditions be satisfied:
\begin{itemize}  
\item[\emph{\textbf{(C1)}}] $R$ is an $(2nr,\eps_1/\sqrt{r},\delta_1)$-JLT. 
\item[\emph{\textbf{(C2)}}] $R$ is a $(r, \eps_2, \delta_2)$-restricted isometry.
\end{itemize}
Then with probability at least $1 - \delta_1 - \delta_2$,
\begin{equation}\label{eq:mybound_bias_cls}
\nnorm{(I - P_{X_R}) X}_F^2 \leq \left( 1 +  \frac{\eps_1^2}{(1 - \eps_2)^4} \right)\nnorm{\Delta_r}_F^2.
\end{equation}
Conditional on the event \eqref{eq:mybound_bias_cls}, the excess risk of CLS \eqref{eq:excess_cls} can be bounded as
\begin{equation}\label{eq:mybound_excess_cls}
\mc{E}(R) \leq \left( 1 +  \frac{\eps_1^2}{(1 - \eps_2)^4} \right) \nnorm{w^*}_2^2 \frac{\nnorm{\Delta_r}_F^2}{n} + \sigma^2 \frac{k}{n}. 
\end{equation}
\end{theo}
\noindent A meaningful interpretation of Theorem \ref{theo:mainresult} and the bound \eqref{eq:mybound_excess_cls} requires an understanding of how
large the number of random projections need to be so that the conditions $\textbf{(C1)}$ and $\textbf{(C2)}$ are satisfied. The next statement addresses this key point.
\begin{prop}\label{prop:subgaussian} Let $R$ have entries drawn i.i.d.~from a zero-mean sub-Gaussian
  distribution and variance $1/k$. If $k = \Omega(\eps_1^{-2} r \{\log(r) + \log(n)\} \vee \eps_2^{-2} \log(\eps_2^{-1}) r\footnote{\label{notertilde} If $r = O(1)$, $r$ can be replaced by $r \vee \log n$ to be make the failure probability small as $n$ grows.})$, then $R^{\T}$ satisfies conditions \emph{(\textbf{C1})}, \emph{(\textbf{C2})} with $\delta_1 = \exp(-c' \log(n \, r))$ and $\delta_2 = \exp(-c \log(\eps_2^{-1}) r\footnoteref{notertilde})$ for absolute constants $c,c' > 0$.
\end{prop}
\noindent According to Proposition \ref{prop:subgaussian}, the bound on the excess risk \eqref{eq:mybound_excess_cls} holds with high probability for $k = \Omega(r \log n)$ sub-Gaussian random projections. The logarithmic factor can potentially be removed, in view of results in Halko et al.~\cite{Halko2011}~for \emph{Gaussian} random
projections. Specifically, in this case Halko et al. show that the expectation of the left hand side of \eqref{eq:mybound_bias_cls} with respect to $R$ can be bounded as  
\begin{equation}\label{eq:froberror_Halko}
\E[\nnorm{(I - P_{X_R}) X}_F^2] \leq \left(1 + \frac{r}{k - r - 1} \right) \nnorm{\Delta_r}_F^2,
\end{equation}
for $k \geq r+2$. In particular, for $k = 2r + 1$, the error in Frobenius norm for approximating the matrix $X$ by $P_{X_R} X$ is within a factor two of the $r$-truncated SVD \eqref{eq:partitioning}.

\subsection{Comparison of PCR and CLS}

\textbf{Statistical performance.} Comparing the bounds on the excess risk \eqref{eq:excessrisk_pcr_bound_dense} and \eqref{eq:mybound_excess_cls} for PCR and CLS, respectively, we see
an agreement in their structure up to constant factors (ignoring the potentially spurious log-factor in the required number of projections
$k$) and the change of $\nnorm{\alpha^*}_{\infty}^2$ to $\nnorm{w^*}_2^2$. As a result, CLS profits from a rapid decay in the sequence of singular
values of $X$ as does PCR, enjoying similarly favorable excess risk bounds under scenarios \textbf{(P)} and \textbf{(E)} given at the end of
\S\ref{sec:pcr}; changes in those bounds are only in terms of constants and the replacement of $\nnorm{\alpha^*}_{\infty}^2$ by $\nnorm{w^*}_2^2$.
The latter may amount to a factor of $\dwn$ in the worst case\footnote{if $n < d$, we may assume without loss of generality that $w^*$ is contained
  in the orthogonal complement of the null space of $X$ so that $\nnorm{\alpha^*}_2^2 = \nnorm{w^*}_2^2$}. However, in light of \eqref{eq:excess_pcr_poly} and
\eqref{eq:excess_pcr_exp} this difference does not have much of an effect as long as the spectrum of $X$ exhibits strong decay.

In spite of this, there is an extreme case in which the ratio of the excess risk of CLS and PCR can be arbitrarily large as can be seen from the exact bound \eqref{eq:excessrisk_pcr_exact}: if $\alpha^*$ happens to be perfectly aligned with the top $r$ singular values so that $\alpha_{r+}^* = 0$, we have $\mc{E}(V_r) = 0$. On the other hand, the column space of $X_R$ does not contain that of $U_r$ unless
$k = \dwn$, hence in this rather specific case CLS falls short of PCR. On the other hand, we are not aware of scenarios in which CLS can substantially improve over PCR.

To be fair, it is worth pointing out that the bound \eqref{eq:mybound_excess_cls} need not always be an improvement over those reviewed in $\S$\ref{sec:cls_existing}, but it yields qualitatively a much better fit if the singular values decay rapidly.       
\vskip2ex
\noindent \textbf{Computational cost.} PCR requires access to the top $r$ left singular vectors of $X$ which is typically done via Krylov subspace methods like Lanczos' algorithm \cite{GolubvanLoan} in $O(n d r)$ flops on average. Assuming that $k = O(r)$, this is comparable to CLS whose dominating operation is given by the matrix-matrix multiplication $X R$. However, as discussed in \cite{Halko2011}, reducing PCR to an $O(n d r)$ operation is problematic as the actual computational complexity can vary significantly depending on subtle spectral properties of $X$: while the nominal complexities of the two approaches are comparable, CLS often requires less runtime. In addition, CLS provides several other advantages. The matrix-matrix multiplication is trivially parellelizable, and can be computed in a single pass over the data. The latter property becomes beneficial once $X$ is too large to fit into the main memory since accessing external memory is slow. Lastly, Theorem \ref{theo:mainresult} and Proposition \ref{prop:subgaussian} are likely extendable to structured JLTs \cite{Ailon2006, AilonLiberty2011, KrahmerWard2011, Tropp2011} for which the cost of forming $X R$ is considerably reduced.    

\section{Further topics}
We here discuss several miscellaneous topics that naturally arise from the analysis of the preceding section. 

\subsection{Assessing the applicability of CLS when being on a computational budget}
\begin{figure}
\begin{center}
  \mbox{\includegraphics[width = 0.39\textwidth]{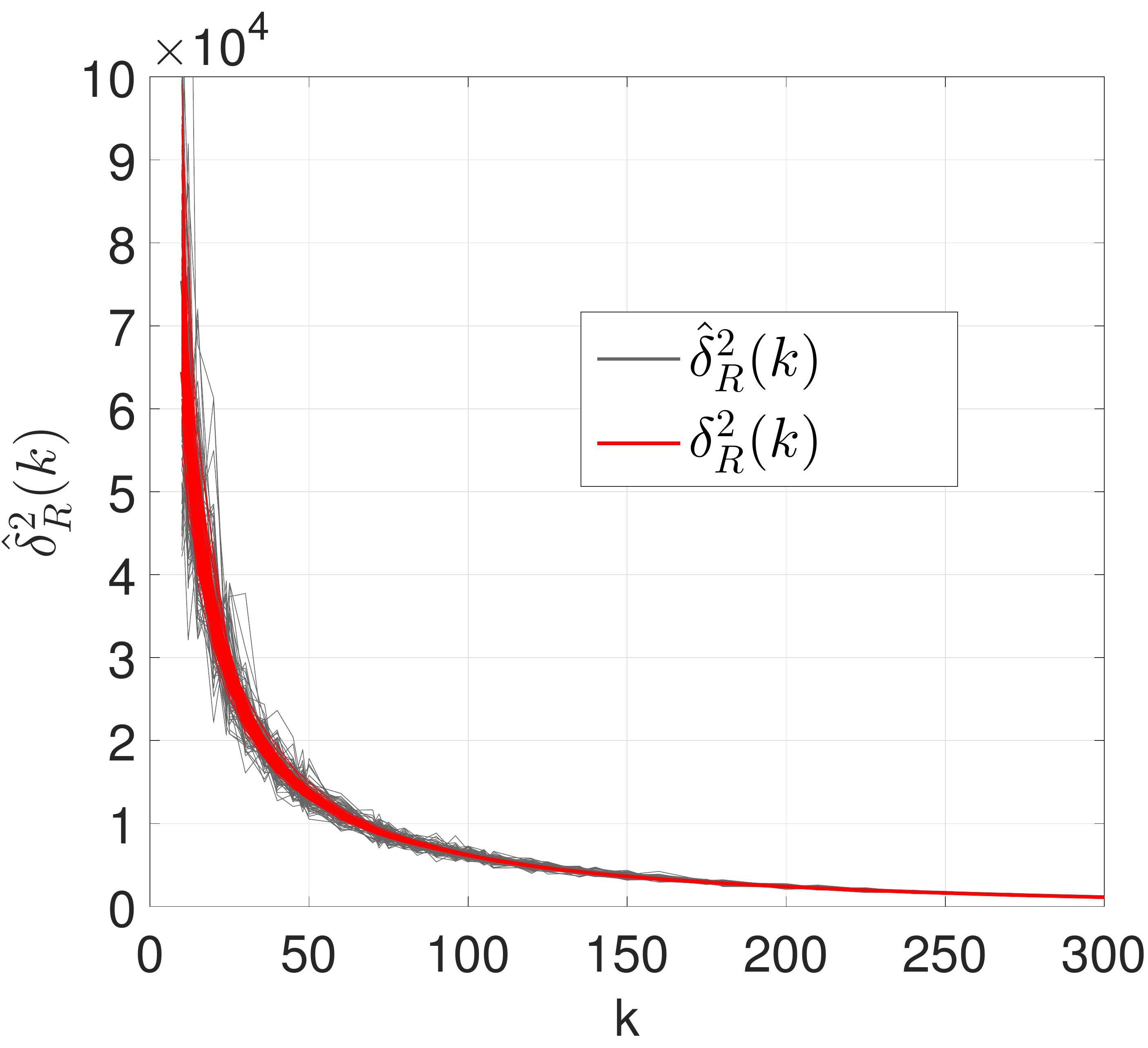} \hskip6ex
    \includegraphics[width = 0.39\textwidth]{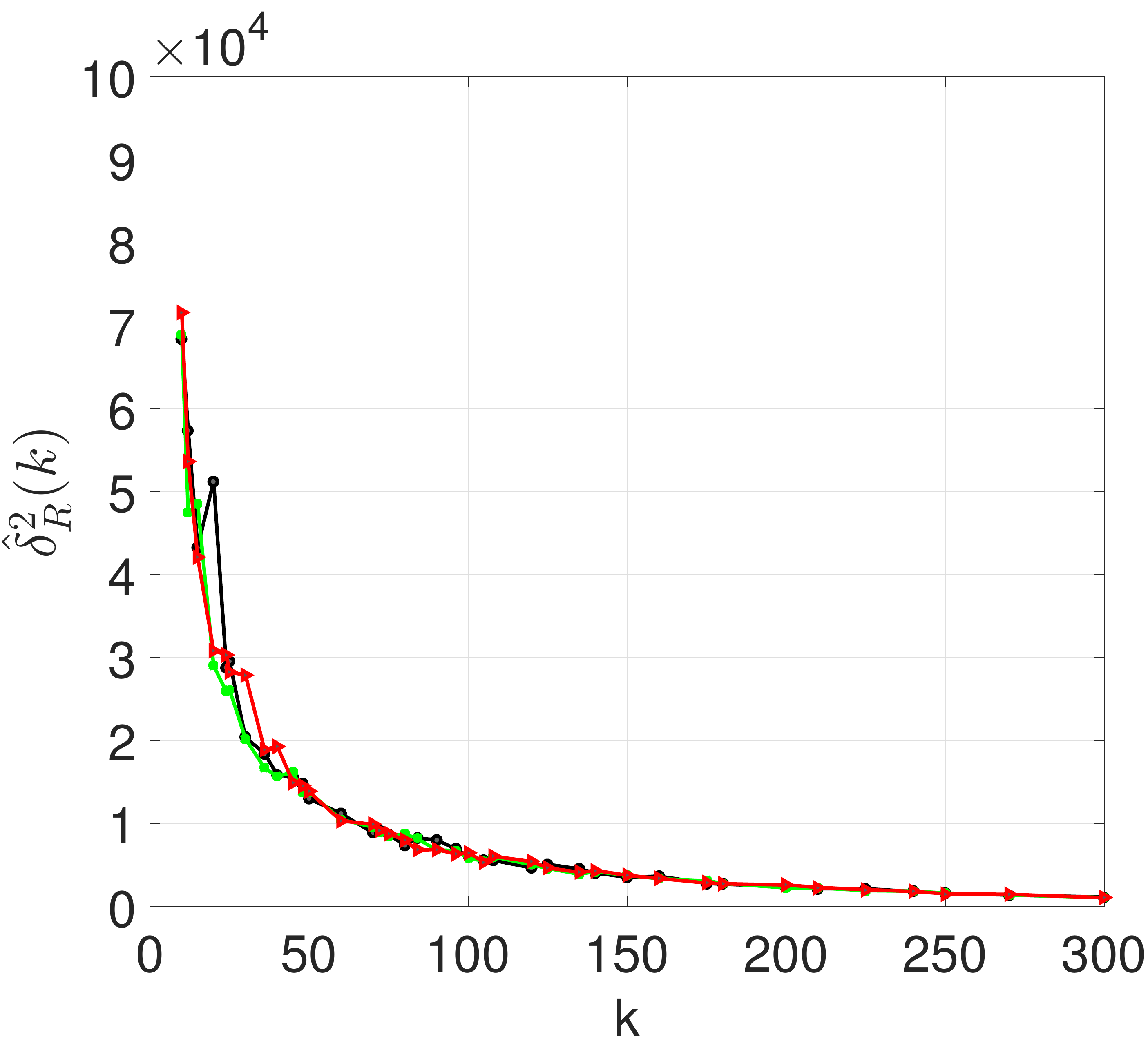}}
\end{center}
\caption{Estimation of $\delta_R^2$ for 100 different realizations of $R$ with i.i.d.~Gaussian entries and a fixed
  matrix $X$ whose singular values decay polynomially (setting (\textbf{P}) in $\S$\ref{sec:pcr}, $q = 2$). Here,
  $\delta_R^2(k)$ is estimated for multiple $k$ simultaneously using the same $\{\omega_l \}_{l = 1}^L$ (here $L = 10$). Left:
  trajectories of $\{ \delta_R^2(k) \}_{k = 1}^{300}$ and $\{ \wwhat{\delta}_R^2(k) \}_{k = 1}^{300}$ for all $100$ realizations of $R$.
  Right: Trajectories of $\{ \wwhat{\delta}_R^2(k) \}_{k = 1}^{300}$ for the first three realizations of $R$.}\label{fig:tailestimation}
\end{figure}

The analysis above indicates that CLS can achieve reasonable statistical performance while using a substantially reduced
number of predictors provided that the singular values of $X$ decay at a fast rate. Verifying such decay seems to require
the SVD though, which would eliminate potential computational advantages of CLS over PCR. Direct evaluation of the quantity $\delta_R^2 = \nnorm{(I - P_{X_R}) X}_F^2$ in \eqref{eq:mybound_bias_cls} is computationally demanding as well: finding the left
singular vectors $\mc{U} \in \R^{n \times k}$ of $X_R$ can be done in  $O(nk^2)$ flops; however, forming $P_{X_R} X = \mc{U} \mc{U}^{\T} X$ requires $O(ndk)$ flops which is of the same order of magnitude as computing $X_R$, the dominating operation in CLS. Hence, we would like to circumvent this operation. We here suggest to recycle the method of random projections in order to get an accurate estimate of $\delta_R^2$ while achieving a reduction to $O(nd)$ flops. The basic idea is to apply $P_{X_R}$ to a small number $L$ of
random elements from the range of $X$ rather than to all its columns.
\begin{prop}\label{prop:randomized_estimation} Consider a collection of $L$ i.i.d.~$d$-dimensional standard Gaussian random vectors $\{ \omega_l \}_{l = 1}^L$ independent
of $R$. Conditional on $R$, consider the following estimator of $\delta_R^2 = \nnorm{(I - P_{X_R}) X}_F^2$:
\begin{equation*}
\wh{\delta}_R^2 = \frac{1}{L} \sum_{l = 1}^L \nnorm{X \omega_l - P_{X_R} X \omega_l}_2^2. 
\end{equation*}
Then, for any $c \in (0,1)$ and any $C > 1$, as long as 
\begin{equation*}
L \geq \max\left\{\frac{16}{(1 - c)^2},  \frac{144}{(C - 1)^2}\right\}
\end{equation*}
it holds that $\p \left( c \delta_R^2 \leq \wh{\delta}_R^2 \leq  C \delta_R^2 \right) \geq 0.96$, where the probability is w.r.t.~$\{ \omega_l \}_{l = 1}^L$ and conditional on $R$. 
\end{prop}
\noindent For example, setting $C = 3$, $c = 1/3$, we would need $L = 36$ to estimate $\delta_R^2$ within a multiplicative factor of $3$ with probability near $1$. Note that computing $P_{X_R} X \omega_l = \mc{U} (\mc{U}^{\T} (X \omega_l))$ for a single $l$ only amounts to $O(nd)$ flops. The constants in Proposition \ref{prop:randomized_estimation} may not necessarily be optimal. In the example of Figure \ref{fig:tailestimation}, we use $L = 10$ random vectors to estimate $\wh{\delta}_R^2(k)$ simultaneously for $1 \leq k \leq 300$. 

\subsection{Column Subsampling vs.~Dense Random Projections}
We can think of CLS as a scheme that picks $k$ random elements from the subspace spanned by
the columns of $X$, and subsequently uses these random elements to predict the response $y$.
The analysis of the previous section asserts that if $k$ is in proper relation to $r$, then
we will do roughly as good as when using the top $r$ principal components as predictors. For this result to hold, Theorem \ref{theo:mainresult} imposes certain restrictions on the matrix $R$, or equivalently, the way the random elements are generated. It is natural to ask whether one can be more flexible in this regard. The simplest approach one can possibly think of would be to select columns from $X$ uniformly at random without replacement, i.e.,
\begin{equation}\label{eq:columnsampling_wo_replacement}
  R = [e_{i_1} \ldots e_{i_k}], \quad i_1 \sim \text{Unif}([d]), \quad i_{\ell} \sim \text{Unif}([d] \setminus \{i_1,\ldots,i_{\ell-1} \}), \; \; \ell=2,\ldots,k,
\end{equation}
where $\{e_1, \ldots, e_d \}$ are the canonical basis vectors of $\R^d$, and $\text{Unif}(\ldots)$ denotes the uniform probability distribution, so that $XR = [X_{i_1} \ldots X_{i_k}]$ is a random column submatrix of $X$. Note that a random matrix $R$ generated according
to \eqref{eq:columnsampling_wo_replacement} fails to be a $(1, \eps, \delta)$-JLT for any $\eps \in (0,1)$, any $\delta \in (0,1/2)$ and any $k \leq \frac{d}{2}$\footnote{To see this, consider the first canonical basis vector $e_1$. Then $\p(R^{\T} e_1 = 0) \geq 1/2$ for all $k \leq d/2$.}, hence the framework used for Theorem \ref{theo:mainresult} does not yield useful results. In general, the excess risk of
column subsampling can be considerably worse than when using a suitable JLT. Both approaches are equivalent in terms of excess risk if i) $n \geq d$ and the spectrum is flat  (scenario \textbf{(F)} in \S\ref{sec:pcr}), or ii) $X$ is a random Gaussian matrix with i.i.d.~$N(0,1)$ entries.
\begin{prop}\label{prop:columsampling_vs_rp} Let $S$ by a random $d \times k$ matrix generated according to \eqref{eq:columnsampling_wo_replacement} and let
$R$ be an $d \times k$ random matrix with i.i.d.~$N(0,1/k)$ entries.
\begin{itemize}
\item[i)] If $n \geq d$ and $X$ has flat spectrum, i.e., $\sigma_1 = \ldots = \sigma_d = \sqrt{n}$, then
      \begin{equation*}
        \E[\mc{E}(R)] = \E[\mc{E}(S)] = \left(1 -  \frac{k}{d} \right) \nnorm{\alpha^*}_2^2 + \frac{k}{n} \sigma^2. 
      \end{equation*}
    \item[ii)] If $n < d$, there exists an $X$ with flat spectrum, i.e., $\sigma_1 = \ldots = \sigma_n = \sqrt{d}$ such that
      \begin{equation*}
        \E[\mc{E}(R)] = \nnorm{\alpha^*}_2^2 \frac{d}{n} \left(1 -  \frac{k}{n} \right) \nnorm{\alpha^*}_2^2  +  \frac{k}{n} \sigma^2 <
        \E[\mc{E}(S)] = \nnorm{\alpha^*}_2^2 \frac{d}{n} \left(1 -  \frac{k}{d} \right) + \frac{k}{n} \sigma^2. 
      \end{equation*}
    \item[iii)] If $X$ is a random $n \times d$ matrix with $N(0,1)$ entries, then $\E[\mc{E}(R)] = \E[\mc{E}(S)]$, where the 
      expectations are w.r.t.~the randomness of both $X$ and $R$ respectively $X$ and $S$. 
\end{itemize}
\end{prop}
\noindent We note that property iii) crucially relies on Gaussianity of the entries of $X$. Proposition \ref{prop:columsampling_vs_rp} will be complemented with numerical results presented in $\S$\ref{sec:experiments} below.   
\subsection{Averaging}\label{sec:averaging}
An idea in the spirit of bagging \cite{Breiman1996} considered in \cite{Thanei2017} is to generate an ensemble of i.i.d.~random projections $\{ R_b \}_{b = 1}^B$ and then average the resulting predictions
\begin{equation*}
\frac{1}{B} \sum_{b = 1}^B X R_b \wh{w}_{R_b} = \frac{1}{B} \sum_{b = 1}^B P_{XR_b}, 
\end{equation*}
where $\wh{w}_{R_b}$ is the least squares estimator given response $y$ and design matrix $X_{R_b}$, $b \in [B]$,
and $\{ P_{XR_b} \}_{b = 1}^B$ are the projections on the respective column spaces of $\{ X R_b \}_{b = 1}^B$.
As $B \rightarrow \infty$, the average of projectors can be expected to behave similarly to $\mc{P}_k = \E[P_{XR}]$. The next
proposition summarizes basic properties of averaging and the operator $\mc{P}_k$.
\begin{prop}\label{prop:averaging} Let $R$ and $\{ R_b \}_{b = 1}^B$ be i.i.d.~random projections. 
\begin{itemize}  
\item[i)] Reduction in bias: $\forall B \geq 1$, $\E[\nnorm{(I - \frac{1}{B}\sum_{b = 1}^B P_{XR_b}) Xw^*}_2^2/n] \leq \E[\nnorm{(I -  P_{XR}) Xw^*}_2^2/n]$.
\item[ii)] Let $R$ have i.i.d.~$N(0,1)$ entries and let $\Delta(k) \coloneq \{z \in \R^{\dwn}: \, \sum_{j = 1}^{\dwn} z_j = k, \, z_j \geq 0, j \in [\dwn] \}$. Then there exists $\{ \eta_j \}_{j = 1}^{\dwn} \in \Delta(k)$ depending only on $\{ \sigma_j \}_{j = 1}^{\dwn}$ s.t.{\small
            \begin{equation*}
            \frac{1}{n} \nnorm{(I  - \mc{P}_k) X w^*}_2^2 =  \frac{1}{n} \sum_{j = 1}^{\dwn} \{ \alpha_j^* \}^2  \sigma_j^2 (1 - \eta_j)^2
              \leq     \frac{1}{n} \E[\nnorm{(I  - P_{XR}) X w^*}_2^2] =  \frac{1}{n} \sum_{j = 1}^{\dwn} \{ \alpha_j^* \}^2  \sigma_j^2 (1 - \eta_j). 
           \end{equation*}}  
\item[iii)] Reduction in variance: $\forall B \geq 1$, $\E\left[\nnorm{\frac{1}{B}\sum_{b = 1}^B P_{XR_b} \xi}_2^2/n \right] \leq \sigma^2 k/n$.
\item[iv)] In the setting of ii), we additionally have
           \begin{equation*}
           \frac{1}{n} \E[\nnorm{\mc{P}_k \xi}_2^2] = \frac{\sigma^2}{n} \tr(\mc{P}_k^2) = \frac{\sigma^2}{n} \sum_{j = 1}^{\dwn} \eta_j^2 = \E\left[\sum_{\ell = 1}^k \cos^2 \theta_{\ell}(\text{\emph{range}}(XR), \text{\emph{range}}(XR')) \right],
         \end{equation*}
         where the last expectation is w.r.t.~$R$ and $R'$ drawn independently, and $\{ \theta_{\ell}(\mc{L}, \mc{L}') \}_{\ell = 1}^k$ denote the canonical angles between two $k$-dimensional subspaces $\mc{L}$ and $\mc{L}'$ \cite{GolubvanLoan}. 
         \end{itemize}
\end{prop}
\noindent Parts ii) and iv) already appear in \cite{Thanei2017} for the $n \geq d$ case with slight differences in presentation. Parts ii) and iv) provide a rough quantification
of the reduction in bias and variance by averaging in the limit $B \rightarrow \infty$ and Gaussian random projections. Since the relationship between the singular values
$\{ \sigma_j \}_{j = 1}^{\dwn}$ and the coefficients $\{ \eta_{j} \}_{j = 1}^{\dwn}$ is not well understood, a more precise connection is yet to be made. In light of Theorem \ref{theo:mainresult},
we know that when choosing $k = \Omega(r \log n)$ large enough, the bias essentially depends only on the tail $\{ \sigma_j^2 \}_{j > r}$. Accordingly, the $\{ \eta_j \}_{j = 1}^r$ must be close
to one, and accordingly the variance term in iv) is at least about $\sigma^2 r /n$. For a choice of $k$ that makes the tail $\{ \sigma_j^2 \}_{j > r}$ (and hence the bias) negligible, averaging
thus does not yield significant benefits. However, when using averaging, the optimal choice of $k$ is guaranteed to be lower than when using a single random projection.

For fixed $k$, the maximum
possible reduction in variance is seen to be a factor $k/(d \wedge n)$, from $\sigma^2 k/n$ to $\sigma^2 k^2 / \{(\dwn) \cdot n \}$, which follows immediately from the representation of the variance
in terms of the $\{ \eta_j \}_{j = 1}^{\dwn}$ in iv) and the fact that the minimum $\ell_2$-norm over the simplex $\Delta(k)$ is attained at its barycenter. It is not hard to see that this maximum reduction
is attained precisely when the spectrum of $X$ is flat. This case remains of limited interest though as the bias is of the same order as without averaging.

The last identity in iv) provides a geometric interpretation of the variance after averaging in terms of the expected sum of squared cosines of the principal angles between to independently generated subspaces $\text{range}(XR)$ and $\text{range}(XR')$. Again, we see that the bias is low, i.e., if $\text{range}(XR)$ well approximates $\text{range}(X)$, two randomly chosen subspaces will be essentially aligned so that the squared cosines evaluate all about one, which implies that averaging will not achieve any substantial reduction in variance. 

\section{Experiments}\label{sec:experiments}
We present the results of experiments with synthetic and real data in order to illustrate and support the main results of the paper, as well as pointing
to some open questions.
\subsection{Synthetic data}\label{subsec:experiments_synthetic}
We start by generating a random $n$-by-$d$ matrix $X_0$ with $n = 1000$, $d = 500$,
where the entries of $X_0$ are drawn i.i.d.~from the $N(0,1)$ distribution. The
SVD of $X_0$ is given by
\begin{equation}\label{eq:X0}
  X_0 = U_0 \Sigma_0 V_0^{\T}, \quad U_0 \in \R^{n \times d}, \;\,\Sigma_0 \in \R^{d \times d}, \;\,V_0 \in \R^{d \times d}. 
\end{equation}
We then replace $\Sigma_0$ by a diagonal matrix $\Sigma$ whose diagonal elements $\{\sigma_j \}_{j=1}^d$
are chosen in a deterministic fashion according to one of the following regimes:
\begin{alignat*}{2}
  &\text{polynomial}:\; &&\sigma_j \propto j^{-q},  \; q \in \{.5,.75,1,1.5, 2, 4 \},\, j \in [d],\\
  &\text{exponential}:\; &&\sigma_j \propto 0.9^j, \;\, j \in [d], 
\end{alignat*}
where the constant of proportionality is determined by the scaling $\sum_{j = 1}^d \sigma_j^2 = n \cdot d$. We subsequently work with $X = U_0 \Sigma V_0^{\T}$, generating data from the model
\begin{equation}\label{eq:model_syn}
y = Xw^* + \sigma \xi,
\end{equation}
where $w^*$ is drawn uniformly from the unit sphere in $\R^d$, $\sigma \in 2^{p}$, $p \in \{-1,-0.5,\ldots,1\}$, and $\xi$ has i.i.d.~standard Gaussian entries.

Given data $(X, y)$, we then perform PCR
with ten different choices of $r$, using an equi-spaced grid of values depending on the regime according to which $X$ has been generated. For CLS, $R$ is chosen as a standard $d$-by-$k$ Gaussian matrix with
$k = \alpha r$, where the over-sampling factor $\alpha \in \{1,1.2,1.5,2,2.5,3 \}$. We conduct 100 independent replications for each regime. Our main interest is in the bias and the prediction error of PCR and CLS:
\begin{alignat*}{3}
  &\nnorm{(I - P_{U_{r}}) X w^*}_2^2/n  \;\;\; &&\text{vs.} \;\;\;  &&\nnorm{(I - P_{X_R}) X w^*}_2^2/n, \\
  &\nnorm{X w^* - X V_r \wh{w}_{V_r}}_2^2/n  \;\;\; &&\text{vs.} \;\;\;  &&\nnorm{X w^* - X_R \wh{w}_R}_2^2/n,
\end{alignat*}
where $\wh{w}_{V_r}$ and $\wh{w}_R$ denote the least squares estimator for data $(X V_r, y)$ and $(X_R, y)$, respectively. In order to compare CLS and column subsampling, we also generate $R$ as a $d \times k$ column submatrix of the identity chosen uniformly at random, cf.~\eqref{eq:columnsampling_wo_replacement}.   

\begin{figure*}[ht!]
  \begin{tabular}{ccc}
  $\qquad\qquad\quad{\small\textbf{Bias}}$  &  $\qquad \qquad \qquad${\small \textbf{Prediction error (low $\sigma$})}  &$\;\;${\small \textbf{Prediction error (high $\sigma$)}}
  \end{tabular}\\
\subfigure{
\includegraphics[width=0.3\textwidth]{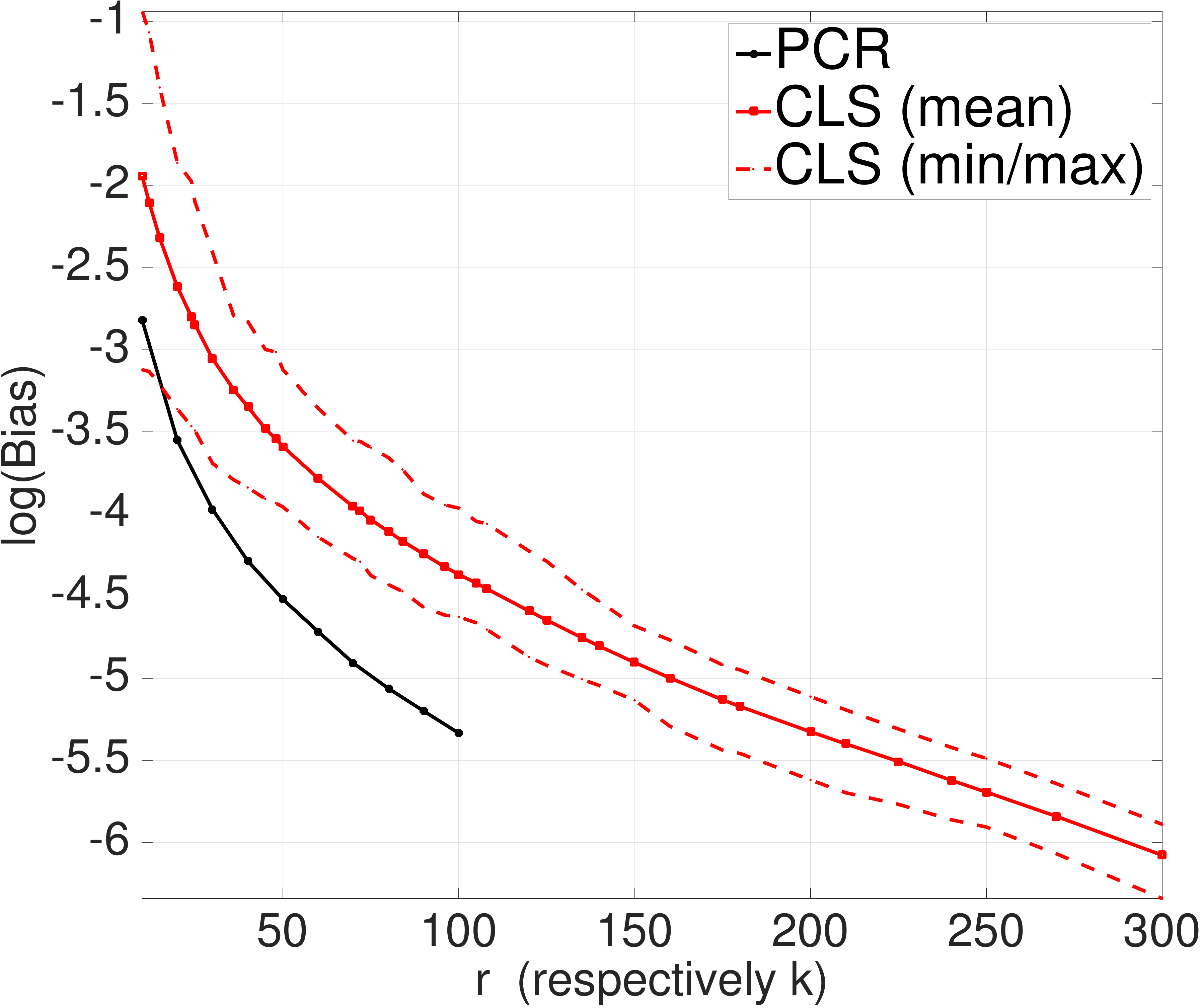}	
\hspace*{0.02\textwidth}	
\includegraphics[width=0.3\textwidth]{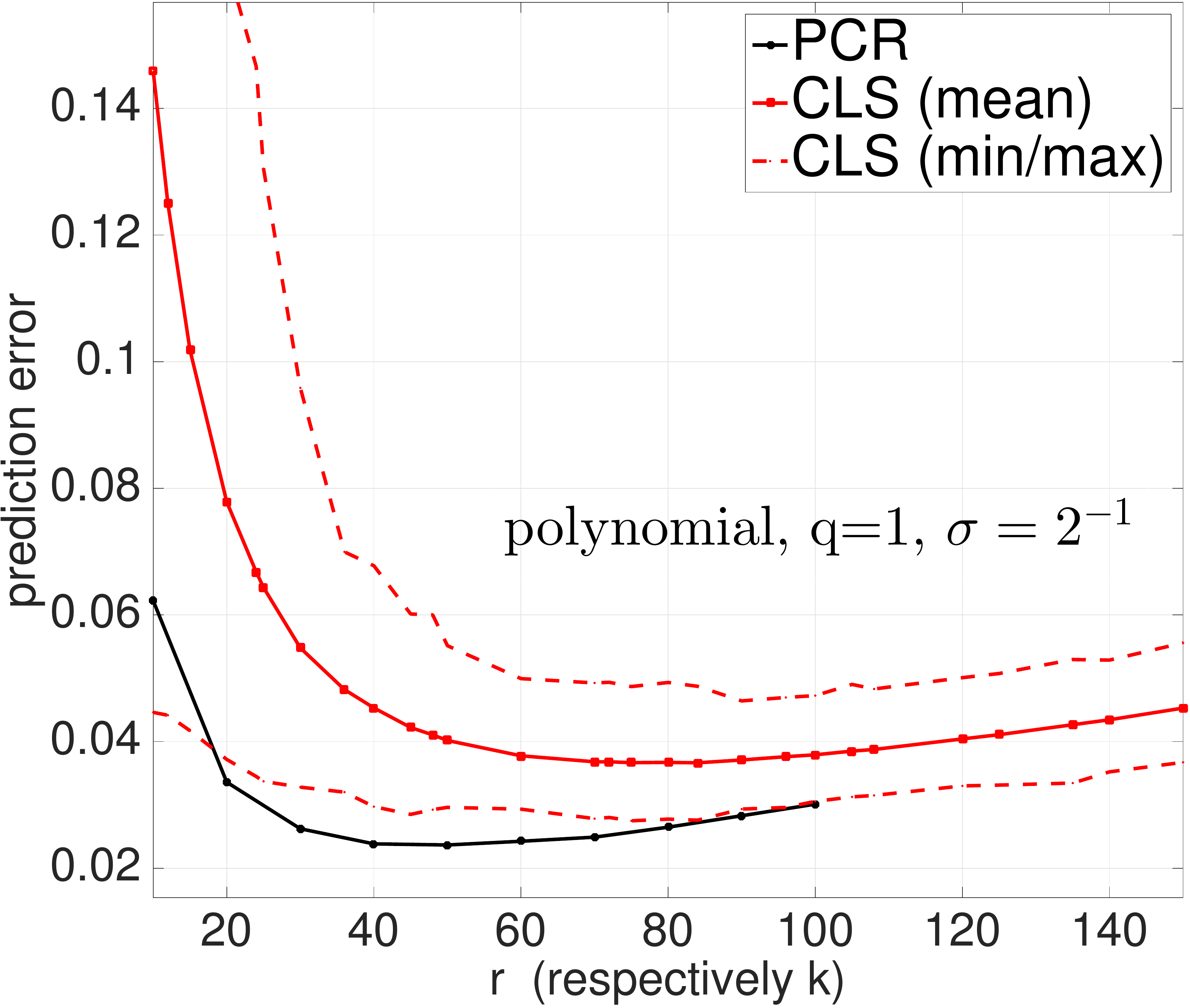}
\hspace*{0.02\textwidth}	
\includegraphics[width=0.3\textwidth]{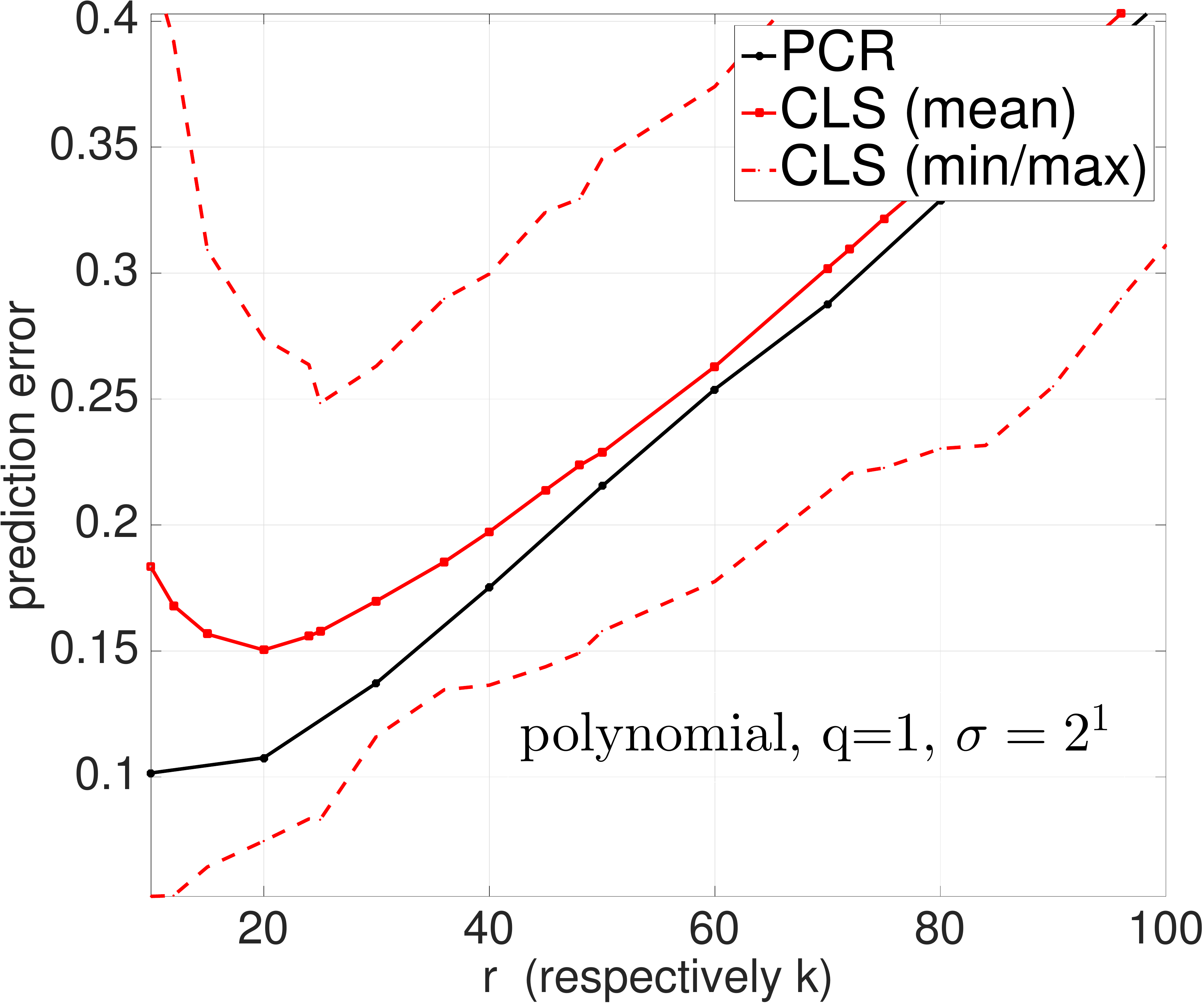}
}
\subfigure{
\includegraphics[width=0.3\textwidth]{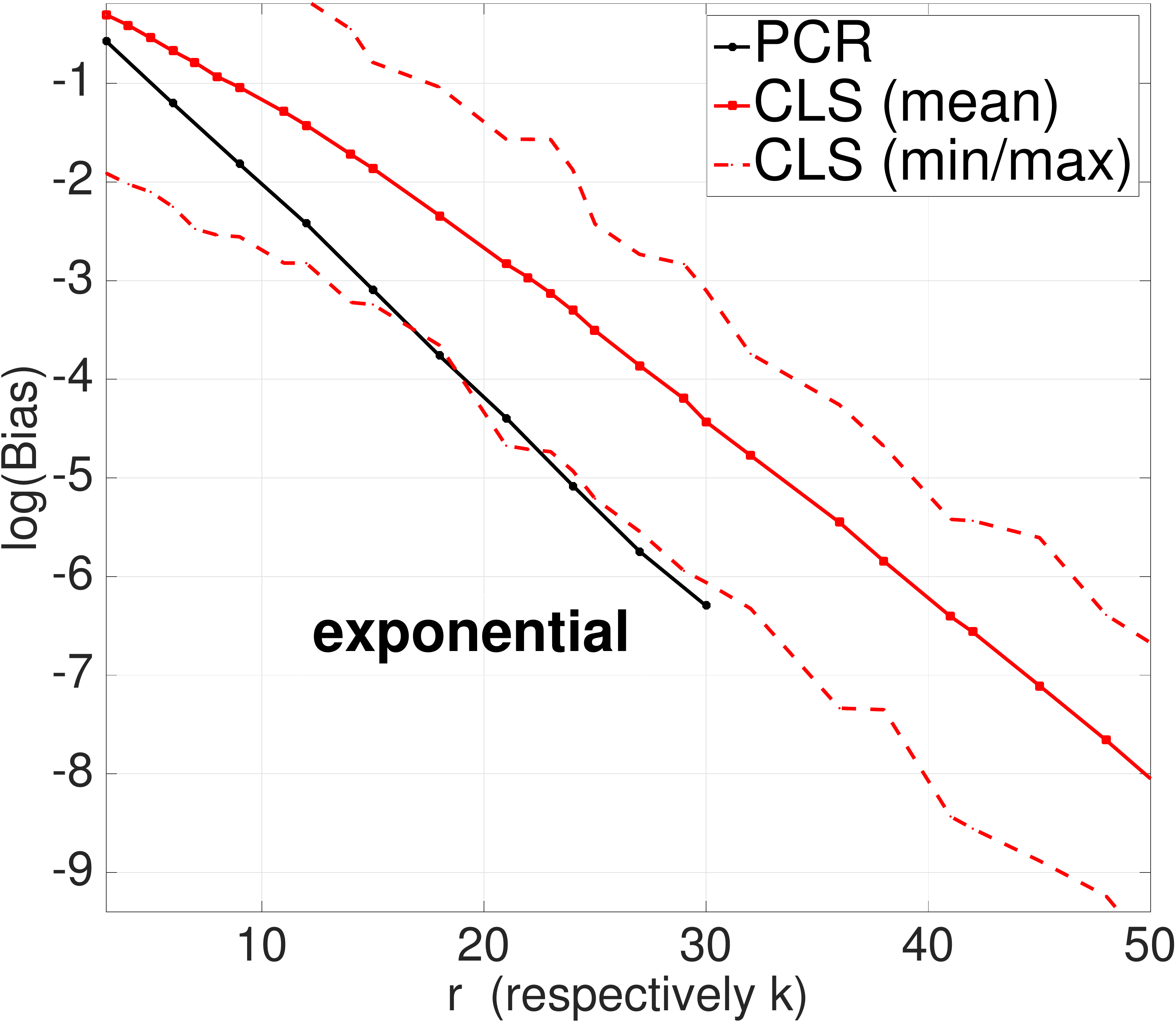}	
\hspace*{0.02\textwidth}	
\includegraphics[width=0.3\textwidth]{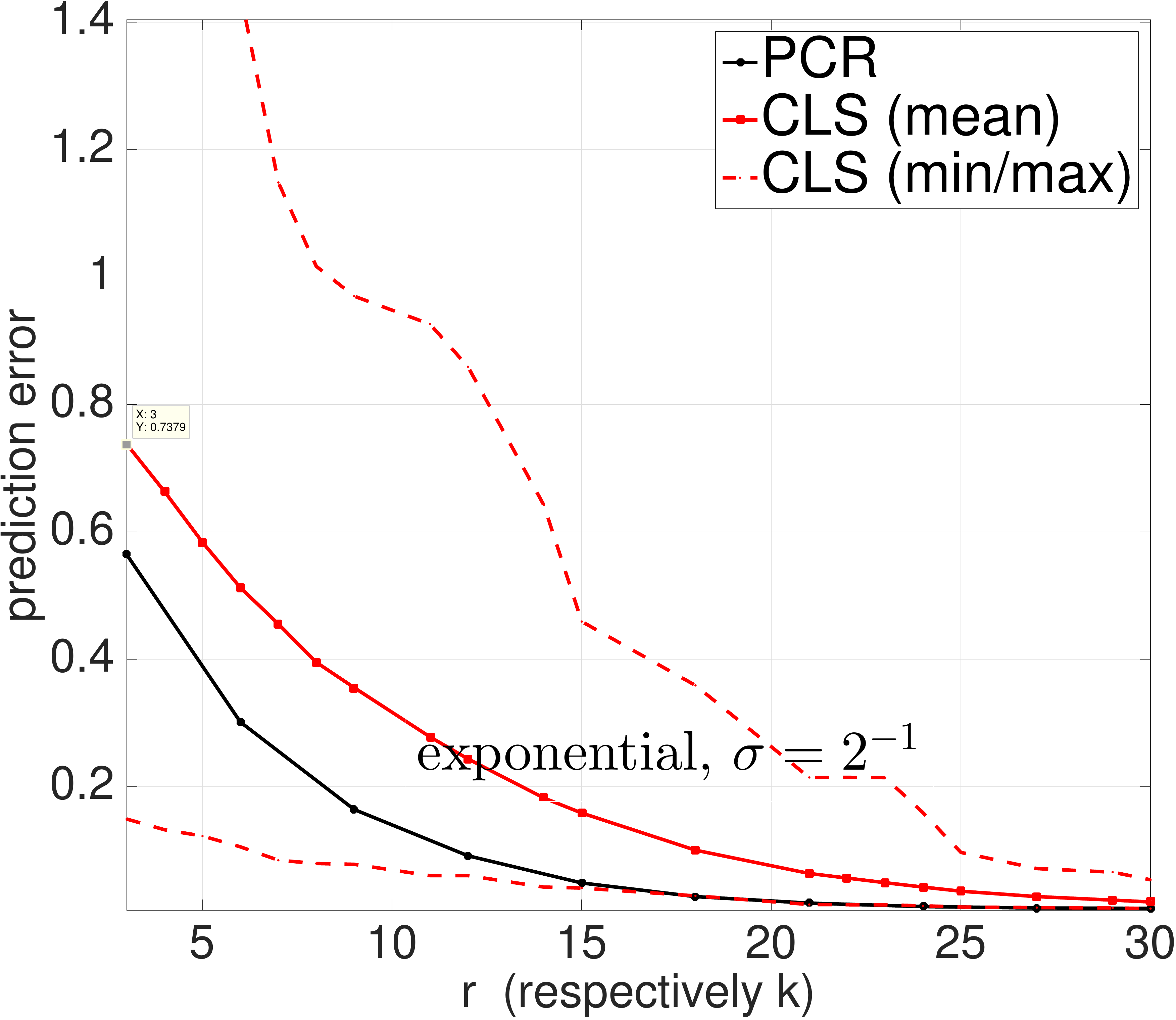}
\hspace*{0.02\textwidth}	
\includegraphics[width=0.3\textwidth]{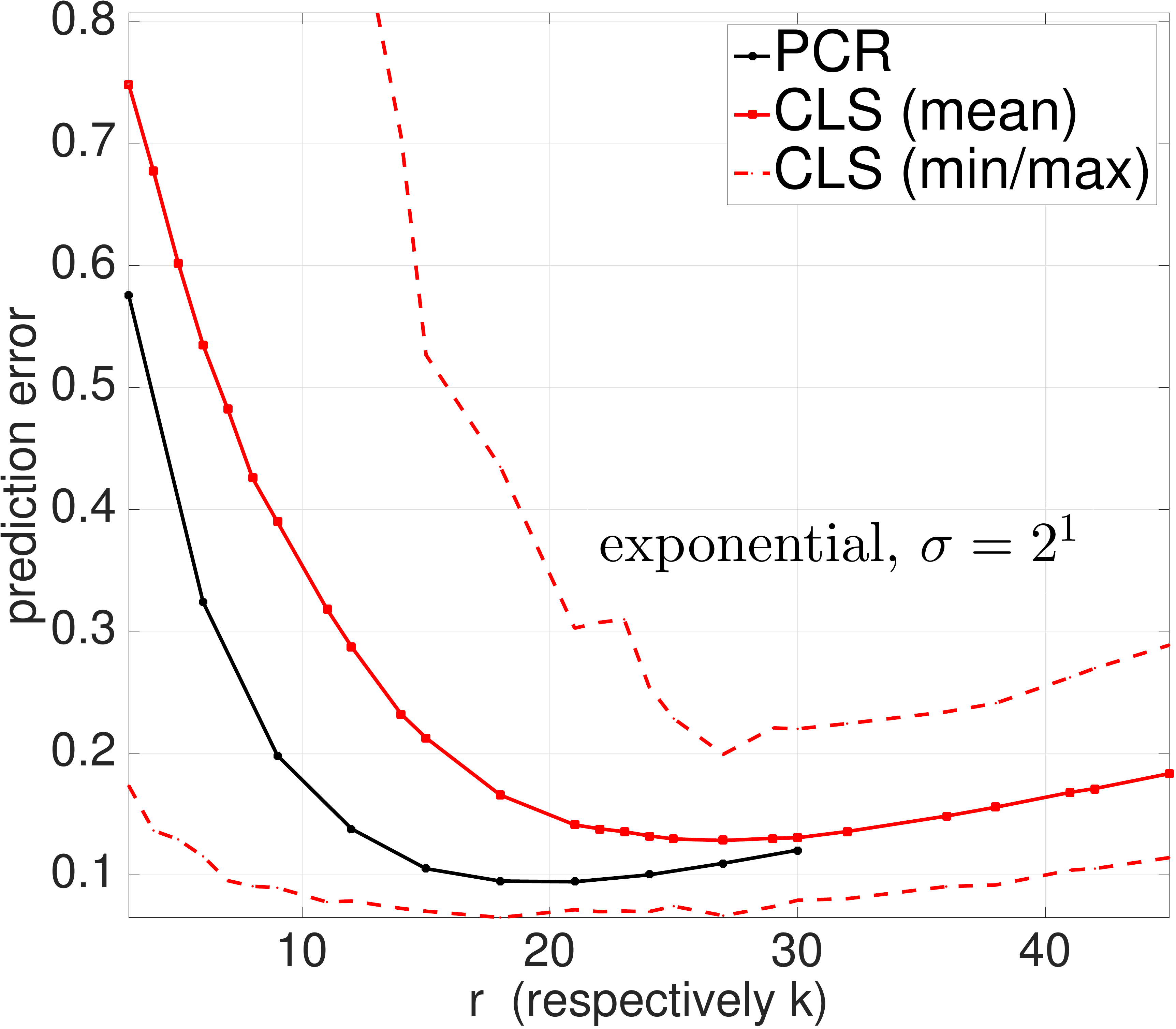}
}
\caption{Results of the synthetic data experiment. From left to right: the bias $\nnorm{(I - P_{X_R}) Xw^*}_2^2 / n$ (log scale), and the mean squared prediction error $\nnorm{Xw^* - X_R \wwhat{w}_R}_2^2/n$ for $\sigma = 1/2$ (middle) and $\sigma = 2$ (right) in dependence of $k$ (horizontal axis) for CLS in relation to PCR. Solid curves are averages, dashed curves minima and maxima (only CLS) over 100 replications. Top: spectrum with polynomial decay ($q = 1$), Bottom: exponential decay. For comparison, the mean squared prediction error of OLS $\sigma^2 d / n$ equals $1/8=.125$ for the middle column and $2$ for the right column.}	
\label{fig:syn}
\end{figure*}

A subset of the results involving two different regimes of decay is shown in Figure \ref{fig:syn}. In the regime of polynomial decay ($q=1$), we observe that the bias of CLS is roughly proportional to that of PCR (or alternatively, we need to choose $k$ as a suitable multiple of $r$ to achieve the same bias). Accordingly, the dip in the prediction error curve occurs for $k = 2r^*$ with $r^* = 40$ yielding the smallest prediction error for PCR. In both low and high noise settings, PCR and CLS improve significantly over OLS in terms of prediction error ($\approx$0.02 and $\approx$0.04 vs.~0.125 and $\approx$ 0.1 and $\approx$ 0.15 vs.~2). In the regime of exponential decay, the bias of CLS is not quite proportional to that of PCR for small values of $r$, but this improves once $r$ reaches 20. Overall, the results agree well with what is suggested by Theorem \ref{theo:mainresult}.

\begin{figure}[ht!]
\begin{center}
\begin{tabular}{cc}
  \includegraphics[width=0.35\textwidth]{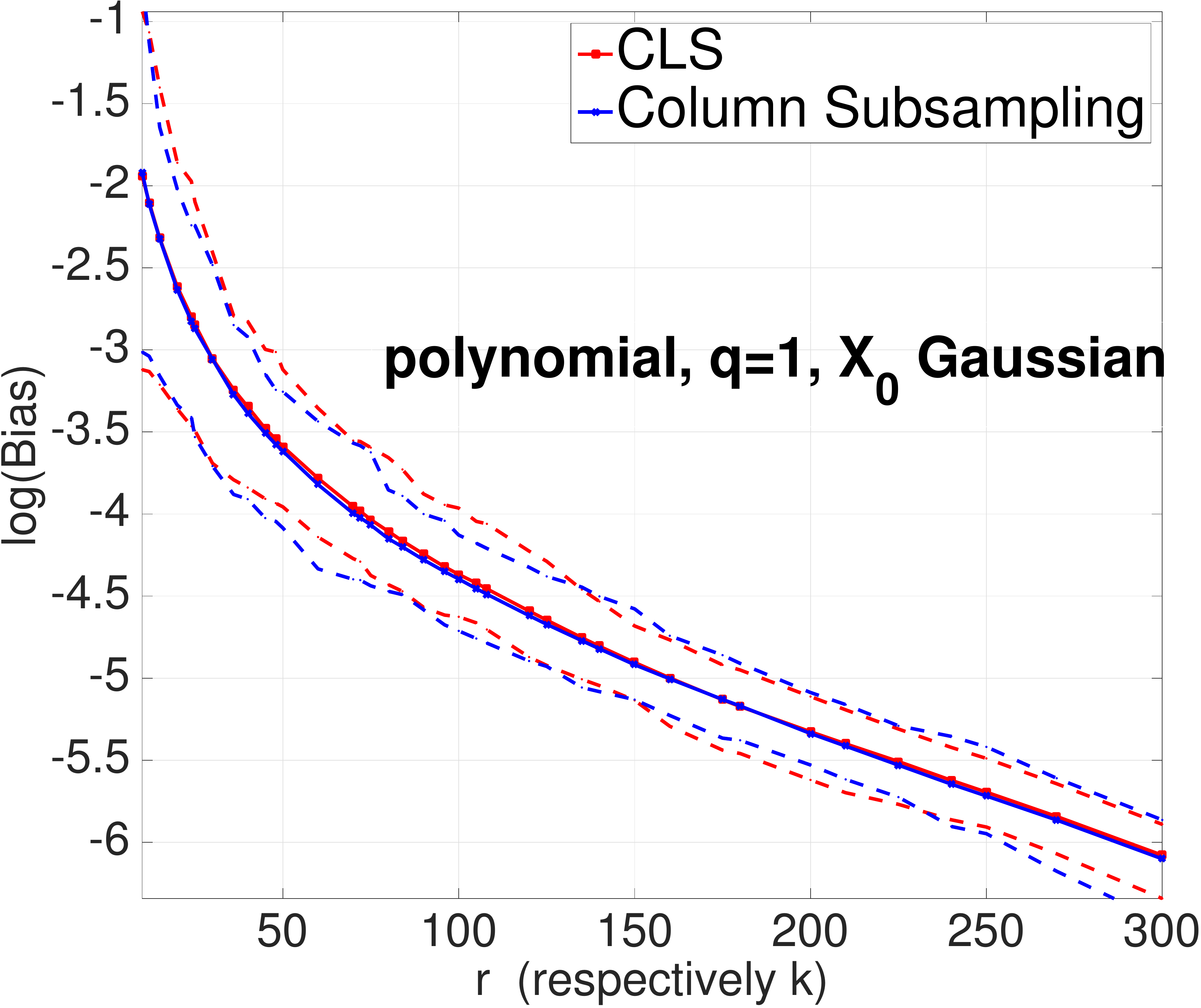}  $\quad$ & $\quad$ 
  \includegraphics[width=0.35\textwidth]{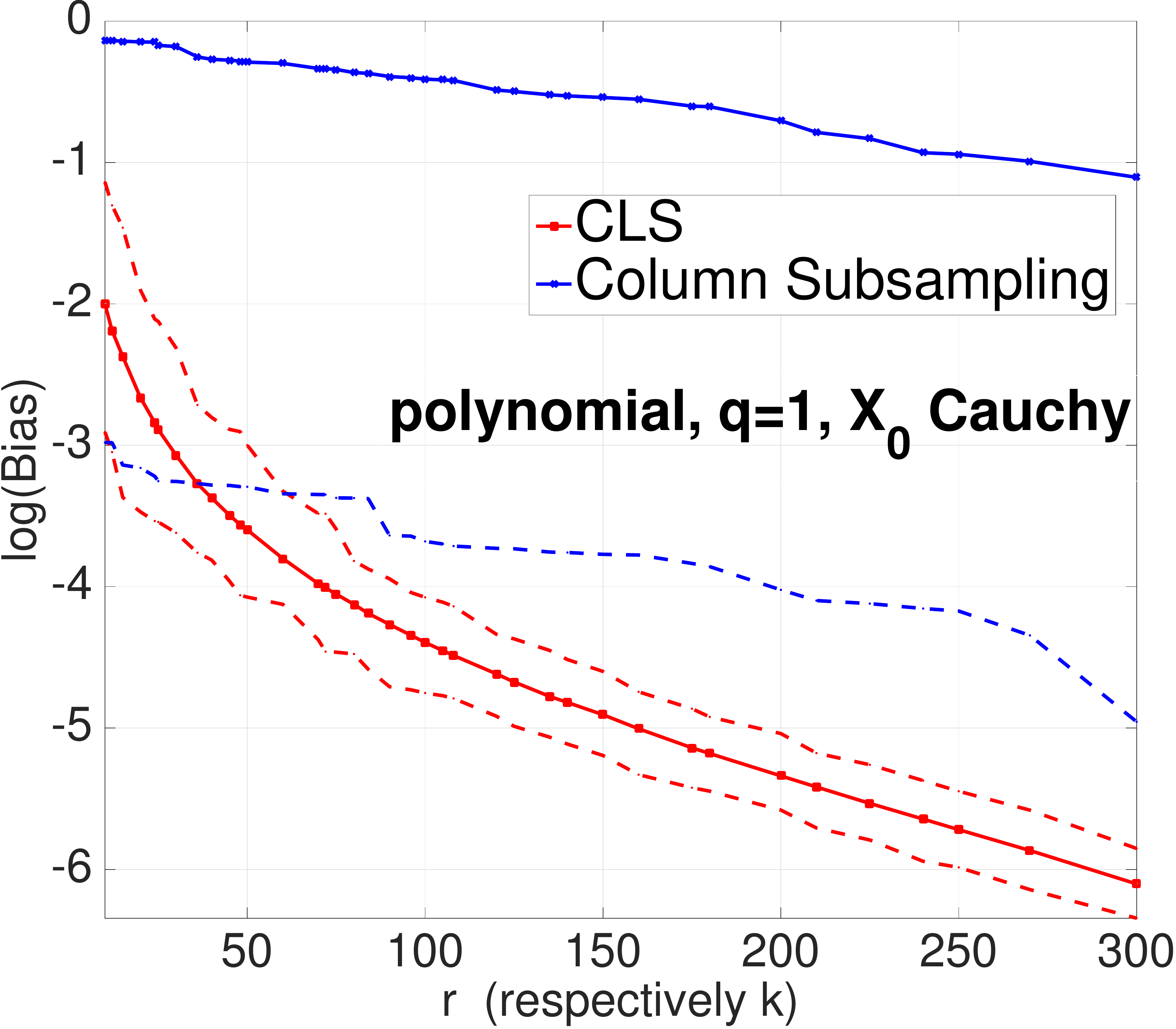}
\end{tabular}
\end{center}
\vspace{-3ex}
\caption{Bias $\nnorm{(I - P_{X_R}) Xw^*}_2^2 / n$ (log scale) when using Gaussian random projections respectively column subsampling. Left: the entries of $X_0$ \eqref{eq:X0} are Gaussian. Right: the entries of $X_0$ \eqref{eq:X0} are Cauchy. Solid curves are averages, dashed curves minima and maxima over 100 replications.}\label{fig:columnsampling_vs_rps}
\end{figure}

Figure \ref{fig:columnsampling_vs_rps} shows that if $X_0$ in \eqref{eq:X0} is generated from a Gaussian distribution, Gaussian random projections and column subsampling perform the same
on average as asserted by Proposition \ref{prop:columsampling_vs_rp}. Interestingly, the performance of column sampling degrades considerably when the entries of $X_0$
are drawn i.i.d.~from the standard Cauchy distribution, whereas Gaussian random projections are not affected by this change.

\subsection{Real data}\label{subsec:experiments_real}

We here consider two data sets from the UCI machine learning repository for illustration purposes. 

\noindent\emph{Twitter social media buzz.} Our data analysis is inspired by that in \cite{Lu2014}. This is a regression problem in which the goal is to predict the popularity of topics as quantified by its mean number of active discussions
given 77 predictor variables such as number of authors contributing to the topic over time, average discussion lengths, number of interactions between authors etc. We here only work with the first $8000$ observations. Several of the original predictor variables as well as the response variable are log-transformed prior to analysis. Following \cite{Lu2014}, we add quadratic interactions which yields $d = 3080$ predictors in total. We consider 50 random partitions into a training set of size $6000$ and a test set of size $2000$ which is used to evaluate the prediction error. The design matrices of the training sets are centered and subsequently scaled to unit norm before performing least squares regression with a centered response and a reduced design matrix obtained from i) a truncated SVD with $r \in \{5,10,\ldots,50,60,\ldots,100,120,\ldots,200 \}$, ii) random Gaussian projections with $k = r \alpha$, where the grid for the factor $\alpha$ is as for the synthetic data, iii) subsets of $k$ columns sampled uniformly at random without replacement. The thus obtained regression coefficients are back-transformed to account for the preliminary centering/scaling step before using them to make predictions on the test set. Approaches i) to iii) are compared in terms of the mean squared prediction error on the test set. For each of the 50 partitions into training and test set, we obtain ten i.i.d.~sets of random projections respectively subsampled columns for ii) respectively iii) and perform regression with each of the resulting reduced matrices, and compute the average as well as the maximum error over each of those ten runs.  
\vskip1ex
\noindent\emph{Blog Feedback.} The task associated with the data set is the prediction of the number comments on blog posts \cite{Buza2014} within a 24-hour time window after a certain base time. The training set consists of $n = 52,397$ observations and originally 280 predictors including meta data about the blogs in which the posts were made, the number of comments received within specific time windows before the base time, number of comments on related posts, bag-of-words data, and the weekday of the post. After eliminating redundant and non-informative predictors, we end up with 114 predictors, several of which are log-transformed. A subset of those are expanded in terms of quadratic and interaction terms which eventually yields $d = 2,589$. The target variable is log-transformed as well. As distinguished from the first case study, this data set comes with a fixed test set of size $7,624$. Evaluation then proceeds as described above, with the only difference that $r \in \{10, 20, 50, 100, \ldots, 500\}$ in order to adjust for a slower rate of decay of the singular values.           
\begin{figure}[h!]
  \centering
  \begin{tabular}{cc}
    {\small \textbf{Twitter}}  &  {\small \textbf{Twitter}} \\
    \includegraphics[height=0.2\textheight]{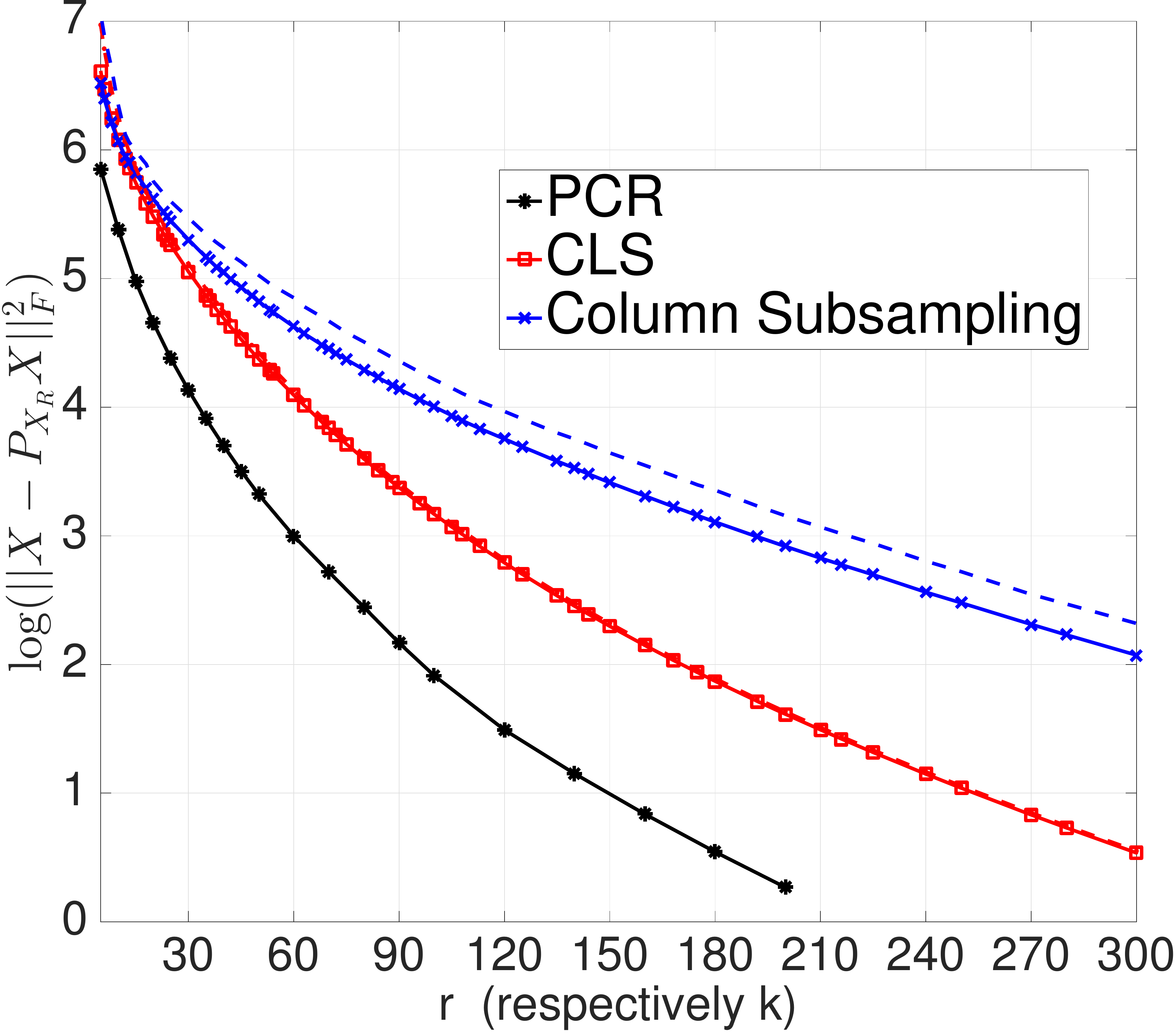} \hspace*{0.02\textwidth} & \hspace*{0.02\textwidth}
                                                                                                                       \includegraphics[height=0.2\textheight]{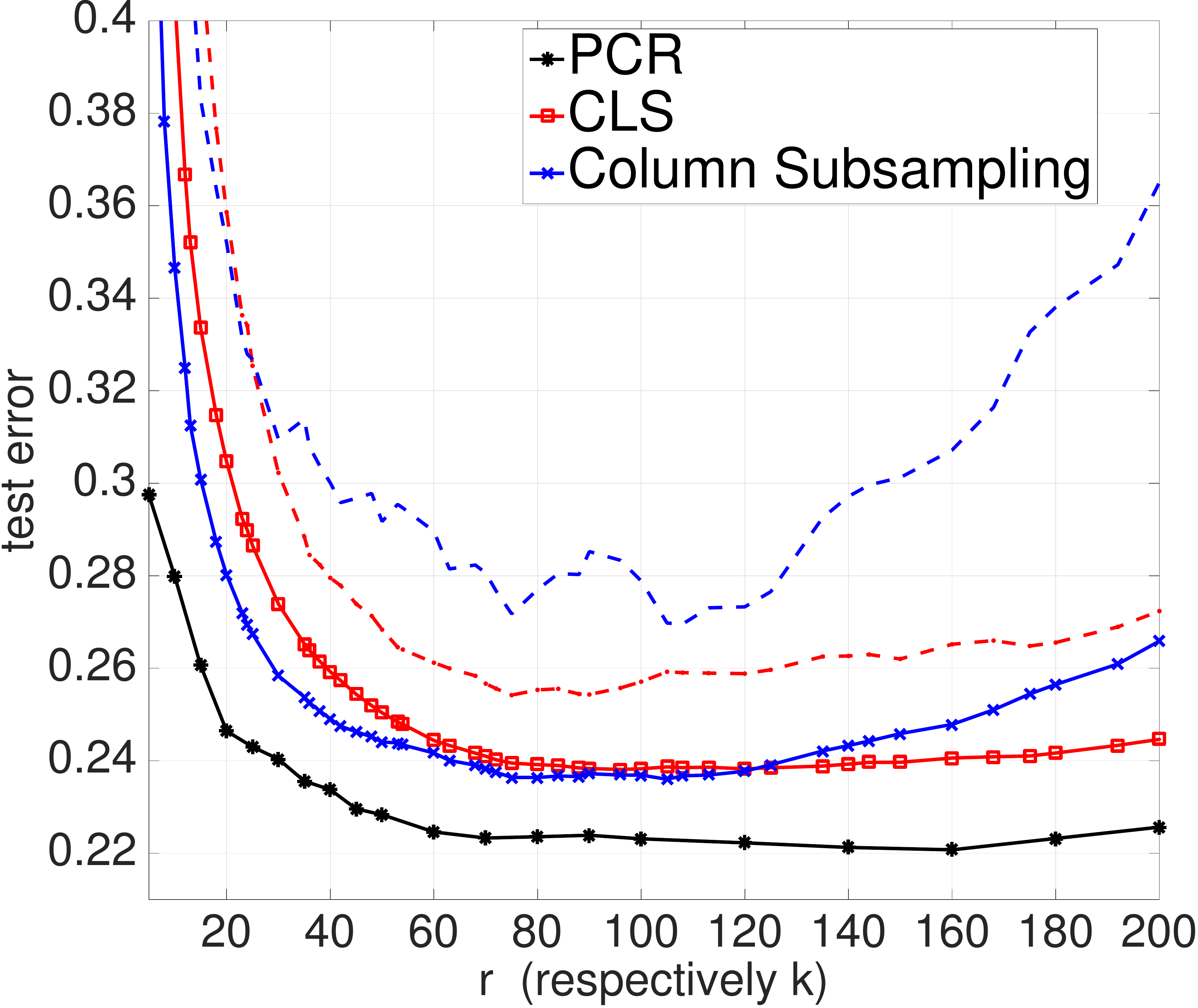} \\
{\small \textbf{Blog Feedback}}  & {\small \textbf{Blog Feedback}} \\    
\includegraphics[height=0.2\textheight]{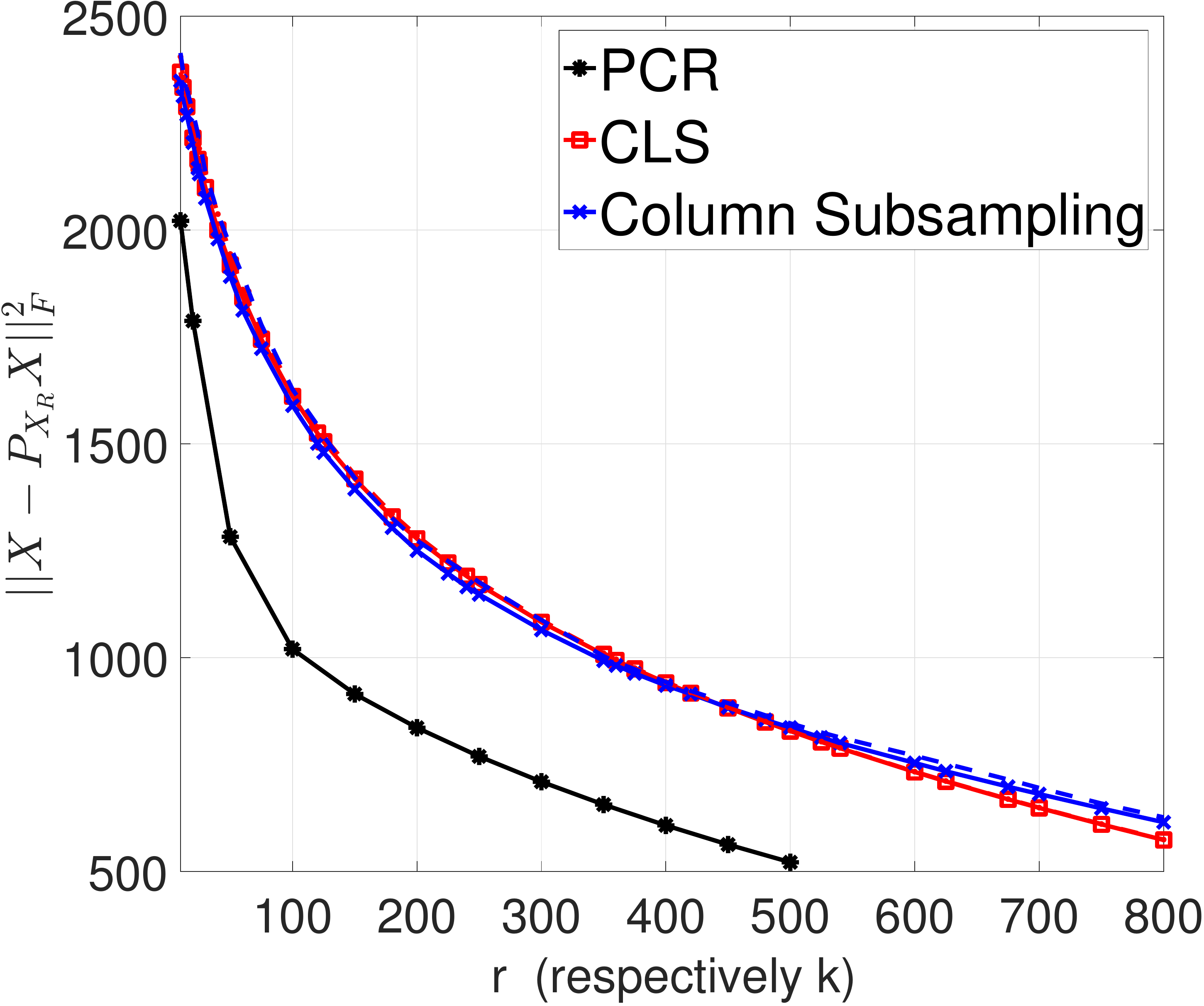} \hspace*{0.02\textwidth} & \hspace*{0.02\textwidth}
                                                                                                                     \includegraphics[height=0.2\textheight]{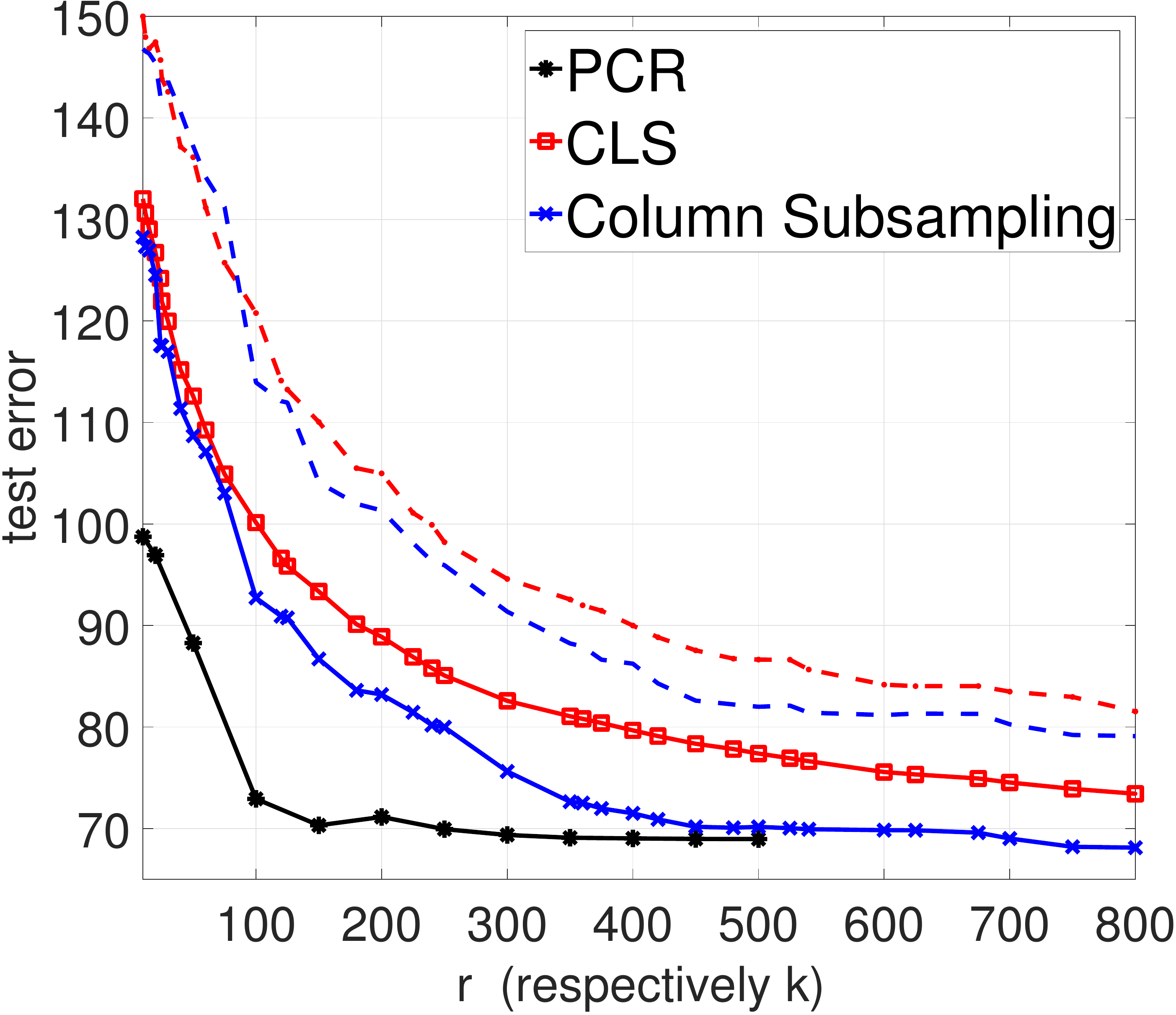}    
  \end{tabular}
  \vspace*{-.7ex}
  \caption{Left: Approximation errors $\nnorm{X - P_{U_r} X}_F^2$ (PCR) respectively $\nnorm{X - P_{X_R} X}_F^2$ (CLS and column subsampling) in dependence of $r$ (or $k$); for ``Twitter'', we use a log-scale and averages over the 50 training/test set partitions. Right: Mean squared prediction errors on the test data vs.~$r$ or $k$. For CLS and column subsampling, we plot both the mean and the maximum over the 10 realizations of $R$.}\label{fig:case_studies}
\end{figure}

The main results are summarized in Figure \ref{fig:case_studies}. The left panels show that the singular values exhibit different rates of decay for the first respectively second data set. For the Twitter data, the decay is not far from linear on a log-scale, whereas the decay is noticeably slower for the Blog Feedback data. As already seen for the synthetic data experiments, CLS requires only a moderate amount of oversampling to achieve the approximation error of PCR $\nnorm{X - P_{U_r} X}_F^2$; we here use this quantity as a surrogate for the bias since $w^*$ is unknown, and in order to establish a connection to result \eqref{eq:mybound_bias_cls} in Theorem \ref{theo:mainresult}. The corresponding quantity of column subsampling is essentially to that of CLS for the Blog Feedback data, but is markedly worse for the Twitter data set. Turning to the right panels, we observe that the test error
of PCR dips for $r = 70$ (Twitter) respectively $r = 150$ (Blog Feedback). The errors of CLS and column subsampling are not far off, but again an increase of $k$ relative to $r$ is indicated to achieve optimal results. Comparing CLS and column sampling, we see that the latter overall performs better particularly for small $k$, but is also less stable in the sense that the gap between the average error and the maximum error (taken over ten different realizations of the random matrix $R$) is substantially larger than for CLS. In particular, for the Twitter data, the maximum error of CLS is uniformly smaller than that of column subsampling. Eventually, the relative performance will depend on properties of $X$ (as illustrated by
Figure \ref{fig:columnsampling_vs_rps}), and likely on possible sparsity of $w^*$ and its interaction with properties of $X$.   
\section{Conclusion}
Regarding the use of random projections in linear regression, the literature has mainly focused on the setting in which the random matrix $R$ is applied
from the left, i.e., $X$ is reduced to $RX$. The converse setting with $X$ being reduced to $XR$ has been studied in several earlier papers as well, however, without establishing a tight link to PCR at the level of achievable excess risk as made herein. Towards the end of our paper, we raise the question how much randomness is needed in generating a random subspace so that such connection holds true. Gaussian random projections induce ``maximum randomness", whereas in column sub-sampling we only consider random subspaces spanned by subsets of the canonical basis vectors. Given the dramatic computational advantages of the latter over ``dense'' random projections it is worth elaborating general conditions under which column sub-sampling can be shown to exhibit similar statistical performance as classical Johnson-Lindenstrauss transforms.

\bibliographystyle{siam}
{
\bibliography{randomizedpcr_long}
}
\appendix

\section{Proof of the lower bound \eqref{eq:bound_thanei_lower}}\label{app:bound_thanei_lower}
Introducing a new set of weights $\{ a_j \}_{j = 1}^{\dwn}$ by the relation $\sigma_j^2 = a_j \tr(\Sigma)$, $\sum_{j = 1}^{\dwn} a_j = 1$,
the $\{ \omega_j \}_{j = 1}^{\dwn}$ \eqref{eq:weights_thanei} can be re-expressed as
\begin{equation*}
\omega_j = \frac{(1 + 1/k) a_j^2 \tr(\Sigma)^2 + (1 + 2/k) a_j \tr(\Sigma)^2 + \tr(\Sigma)^2 / k}{(k + 2 + 1/k) a_j^2 \tr(\Sigma)^2 + 2(1 + 1/k) a_j \tr(\Sigma)^2 + \tr(\Sigma)^2 / k}, \;\; j \in [\dwn].
\end{equation*}
Factoring out the term $\tr(\Sigma)^2$, we arrive at (cf.~Equation (12) in \cite{Thanei2017})
\begin{equation*}
\omega_j = \frac{(1 + 1/k) a_j^2  + (1 + 2/k) a_j  +  1 / k}{(k + 2 + 1/k) a_j^2  + 2(1 + 1/k) a_j  +  1 / k}, \;\; j \in [\dwn]
\end{equation*}
After some algebra, one obtains the following more compact representation:
\begin{align*}
  \omega_j = \frac{1 + a_j}{1 + a_j + a_j k} = \frac{1/a_j + 1}{1/a_j + 1 + k}, \;\; j \in [\dwn]. 
\end{align*}
which is a monotonically decreasing function in $0 \leq a_j \leq 1$, whose minimum is achieved
for $a_j = 1$. Consequently, $\omega_j \geq \frac{2}{2 + k}$, $j \in [\dwn]$, and we conclude \eqref{eq:bound_thanei_lower}.

\section{Proof of Proposition \ref{prop:subgaussian}}
Let $R$ be a random matrix of i.i.d.~sub-Gaussian random variables with zero mean and variance $1/k$.

Regarding \textbf{(C1)}, it is shown in \cite{Matousek2008} (see Theorem 3.1 and the proof therein; cf.~also \cite{IndykNaor2007}) that for any fixed $v \in \R^d$ and $\eps' \in (0,1)$
\begin{equation*}
  \p \left((1-\eps') \nnorm{v}_2^2 \leq \nnorm{R^{\T} v}_2^2 \leq (1+\eps') \nnorm{v}_2^2 \right)
  \leq 2 \exp(-c_0 (\eps')^2 k). 
\end{equation*}
for some absolute constant $c_0 > 0$. It hence follows from the union bound that for any set $\mc{S}$ of vectors in $\R^d$,
$|\mc{S}| = 2 n \cdot r$, 
\begin{equation*}
  \p \left(\forall v \in \mc{S}:\; (1-\eps') \nnorm{v}_2^2 \leq \nnorm{R^{\T} v}_2^2 \leq (1+\eps') \nnorm{v}_2^2 \right)
  \leq \exp(-c_0 (\eps')^2 k + \log(4 n r)).
\end{equation*}
Setting $\eps' = \eps_1 / \sqrt{r}$ for $\eps_1 \in (0,1)$, it follows that for $k = \Omega(\eps_1^{-2} r \log(nr))$, condition \textbf{(C1)} holds with $\delta_1 = \exp(-c' \log(n r))$. 

Turning to \textbf{(C2)}, it follows from arguments in \cite{Baraniuk2006} (cf.~Lemma 5.1 therein) that for any fixed subspace $\mc{V}$ of dimension $r$ in in $\R^d$, $r < k$, 
\begin{equation*}
(1 - \eps_2) \nnorm{v}_2 \leq \nnorm{R^{\T} v}_2 \leq (1 + \eps_2) \nnorm{v}_2 \quad \text{for all $v \in \mc{V}$},
\end{equation*}
with probability at least 
\begin{align*}
1 - 2 (12 / \eps_2)^r \exp(-c_0 \eps_2^2 k) = 1 - \exp \left(-c_0  \eps_2^2 k + r \log(12 / \eps_2)  + \log(2)  \right),
\end{align*}
Hence, for $k = \Omega(\eps_2^{-2} \log(\eps_2^{-1}) r )$, \textbf{(C2)} holds with $\delta_2 = \exp(-c \log(\eps_2^{-1}) r)$. 
This concludes the proof of the proposition.

\section{Proof of Theorem \ref{theo:mainresult}}

Before going into the proof, let us introduce a bit more of notation. Below,  $\mc{T}_s(M)$
denotes the best rank-$s$ approximation of a matrix $M$ with respect to Frobenius norm which
can be obtained from a truncated SVD, cf.~\eqref{eq:excessrisk_pcr_bound_dense}.  Moreover,
we write $M^-$ for the Moore-Penrose pseudoinverse of a matrix $M$. The $j$-th column of $M$ is
denoted by $M_{:,j}$. $\nnorm{M}_2$ denotes the spectral norm. 
\vskip1ex
Note that in view of \eqref{eq:excess_cls} 
\begin{equation*}
\mc{E}(R) = \nnorm{(I - P_{X_R}) Xw^*}_2^2/n  +\sigma^2 \text{rank}(X_R)/n,  
\end{equation*}
we have
\begin{align*}
  \mc{E}(R)  &\leq \left( \nnorm{(I - P_{X_R})}_2^2/n \right) \nnorm{w^*}_2^2 + \sigma^2 k/n \\
                 &\leq \left(\nnorm{(I - P_{X_R})}_F^2/n \right) \nnorm{w^*}_2^2 + \sigma^2 k/n,
\end{align*}
so that \eqref{eq:mybound_excess_cls} immediately follows from \eqref{eq:mybound_bias_cls}. In the
sequel, we hence prove \eqref{eq:mybound_bias_cls}, following the strategy of the proof of Theorem 14 in \cite{Sarlos2006}. 
The proof can be partitioned into three basic steps.
\vskip1ex
\noindent \underline{\emph{\textbf{Step 1.}}}
\begin{lemmaApp}\label{lem:step1_1}  We have
\begin{equation}\label{eq:matrixapprox_1}
\nnorm{X  - P_{X_R} X}_F^2  \leq \nnorm{X  - \mc{T}_r(P_{X_R} X)}_F^2. 
\end{equation}
\end{lemmaApp}
\begin{bew} Observe that according to the definition of $P_{X_R}$, we have  
\begin{equation}\label{eq:matrixapprox_1b}
\min_{B \in \R^{k \times d}} \nnorm{X - X_R B}_F^2 = \nnorm{X - P_{X_R} X}_F^2.  
\end{equation} 
Let $B^* \in \R^{k \times d}$ denote a minimizer of the optimization problem on the l.h.s.~of \eqref{eq:matrixapprox_1b} such
that $P_{X_R} X  = X_R B^*$. Let the SVD of that matrix be given by 
\begin{equation*}
X_R B^* = \underset{n \times d}{\Upsilon} \, \underset{d \times d}{\Xi} \, \underset{d \times d}{\Psi^{\T}}.
\end{equation*}
Denote by $M_r \in \R^{d \times d}$ the diagonal matrix whose first $r$ diagonal entries are equal to one and zero else. Then
$\mc{T}_r(P_{X_R} X) = X_R B^* M_r = X_R \wt{B}$. Since $\wt{B} = B^* M_r$ is a feasible solution for the minimization problem \eqref{eq:matrixapprox_1b}, we conclude \eqref{eq:matrixapprox_1}.   
\end{bew}

\begin{lemmaApp}\label{lem:step1_2} We have
\begin{equation}\label{eq:matrixapprox_2}
\nnorm{X  - \mc{T}_r(P_{X_R} X)}_F^2 \leq \nnorm{X - \Pi X}_F^2,
\end{equation}
where $\Pi$ is the orthogonal projection on the subspace spanned by the columns of $\Phi_r = P_{X_R} \mc{T}_r(X)$, i.e.
\begin{equation}\label{eq:Pi}
\Pi = P_{\Phi_r} = P_{P_{X_R} \mc{T}_r(X)}.  
\end{equation}
\end{lemmaApp}

\begin{bew} Consider the following optimization problem:
\begin{equation*}
\min_{\text{rank}(B) \leq r} \nnorm{X - X_R B}_F^2. 
\end{equation*}
Then any minimizer $B^*$ of the above problem satisfies $X_R B^* = \mc{T}_r(P_{X_R} X)$ (see Proposition 1 and Lemma 14 in \cite{Bunea2011}). Noting that $\Pi = X_R M$ for some matrix $M \in \R^{k \times d}$ with $\text{rank}(M) \leq r$ 
(as $ \mc{T}_r(X)$ has rank no more than $r$), $M$ is feasible for the above optimization problem, and we conclude \eqref{eq:matrixapprox_2}.
\end{bew}

We conclude \emph{\textbf{Step 1.}} by combining Lemmas \ref{lem:step1_1} and \ref{lem:step1_2}:
\begin{equation}
\nnorm{X  - P_{X_R} X}_F^2  \leq \nnorm{X - \Pi X}_F^2
\end{equation}
with $\Pi$ defined in \eqref{eq:Pi}.  
\vskip2ex
\noindent \underline{\emph{\textbf{Step 2.}}}
\vskip1ex
\noindent
In the second step, we decompose $\nnorm{X - \Pi X}_F^2$ into two parts: an ``easy'' part and one more delicate part that requires sophisticated analysis. Recalling \eqref{eq:partitioning}, we have
\begin{align}
\nnorm{X - \Pi X}_F^2 &=  \nnorm{U \Sigma V^{\T} - \Pi U \Sigma V^{\T}}_F^2 \notag\\
                   &=  \nnorm{U \Sigma  - \Pi U \Sigma}_F^2 \notag\\
                      &=  \nnorm{U_{r} \Sigma_{r}  - \Pi U_{r} \Sigma_{r}}_F^2 +\nnorm{U_{r+} \Sigma_{r+}  - \Pi U_{r+} \Sigma_{r+}}_F^2\notag\\
                      &= \nnorm{U_{r} \Sigma_{r}  - \Pi U_{r} \Sigma_{r}}_F^2   + \nnorm{(I - \Pi) U_{r+} \Sigma_{r+}}_F^2  \notag\\  
&\leq \nnorm{U_{r} \Sigma_{r}  - \Pi U_{r} \Sigma_{r}}_F^2 + \nnorm{U_{r+} \Sigma_{r+}}_F^2 \notag\\
&= \underbrace{\nnorm{U_{r} \Sigma_{r}  - \Pi U_{r} \Sigma_{r}}_F^2}_{\text{part requiring special treatment}} + \underbrace{\nnorm{X - \mc{T}_r(X)}_F^2}_{\text{part that we need (up to constant)}} \label{eq:step2_end}
\end{align}
where the inequality follows from the fact that $I - \Pi$ is an orthogonal projection.

\vskip1ex
\noindent\underline{\emph{\textbf{Step 3.}}}
\vskip1ex
\noindent
It remains to bound
\begin{equation*}
\nnorm{U_{r} \Sigma_{r}  - \Pi U_{r} \Sigma_{r}}_F^2 = \nnorm{\mc{T}_r(X) - \Pi \mc{T}_r(X)}_F^2.  
\end{equation*}
Let us write $C^* = X_R^{-}$ and $\wt{C} = (\mc{T}_r(X) R)^-$. Note that for any matrix $M$ of appropriate dimension, we have
\begin{equation}\label{eq:pseudo1}
\nnorm{M - P_{X_R} M}_F^2 = \min_{C \in \R^{k \times n}} \nnorm{M - X_R C}_F^2 =  \nnorm{M - X_R C^*}_F^2 \leq \nnorm{M - X_R \wt{C}}_F^2. 
\end{equation}
Moreover, observe that according to the definition of $\Pi$ in \eqref{eq:Pi}
\begin{equation}\label{eq:doublepi}
\Pi \mc{T}_r(X) = P_{P_{X_R} \mc{T}_r(X)} \mc{T}_r(X) = P_{X_R} \mc{T}_r(X). 
\end{equation}
Using \eqref{eq:pseudo1} and \eqref{eq:doublepi}, we obtain that
\begin{align}
\nnorm{\mc{T}_r(X) - \Pi \mc{T}_r(X)}_F^2 &=  \nnorm{\mc{T}_r(X) - X_R (X_R)^- \mc{T}_r(X)}_F^2 \notag\\ 
                                         &\leq \nnorm{\mc{T}_r(X) - X_R \{\mc{T}_r(X) R\}^{-} \mc{T}_r(X)}_F^2  \notag\\
                                         &=\nnorm{\mc{T}_r(X)^{\T} - \mc{T}_r(X)^{\T} \{ R^{\T} \mc{T}_r(X)^{\T} \}^{-} R^{\T} X^{\T}}_F^2 \label{eq:matrixapprox_3}
\end{align}
Define
\begin{align}\label{eq:leastsquares_A_b}
  b_i = (X^{\T})_{:,i} \in \R^d, \;\;i \in [n], \qquad \text{and} \;\; A = \mc{T}_r(X)^{\T} \in \R^{d \times n},  
\end{align}
and consider the least squares problems
\begin{equation*}
\min_{\lambda_i} \nnorm{b_i - A \lambda_i}_2^2
\end{equation*}
with minimizer $\lambda_i^*$, $i = 1,\ldots,n$, and the corresponding \emph{sketched regression problems} with sketching matrix $R^{\T}$:
\begin{equation*}
\min_{\lambda_i} \nnorm{R^{\T} b_i - R^{\T} A \lambda_i}_2^2, 
\end{equation*}
with minimizer $\wt{\lambda}_i$, $i=1,\ldots,n$. It is straightforward to show that
\begin{equation*}
A \lambda_i^* = (\mc{T}_r(X)^{\T})_{:,i}, \quad i \in [n]. 
\end{equation*}
For the sketched regression problems, an optimal set of coefficients is given by
\begin{equation*}
\wt{\lambda}_i = \{ R^{\T} \mc{T}_r(X)^{\T} \}^{-} R^{\T} (X^{\T})_{:,i}, \quad i \in [n], 
\end{equation*}
so that
\begin{equation*}
A \wt{\lambda}_i = \mc{T}_r(X)^{\T} \{ R^{\T} \mc{T}_r(X)^{\T} \}^{-} R^{\T} (X^{\T})_{:,i}, \quad i \in [n]. 
\end{equation*}
Identifying terms, we see that the right hand side in \eqref{eq:matrixapprox_3} can be written as
\begin{align}\label{eq:normbetai}
\begin{split}
  &\nnorm{\mc{T}_r(X)^{\T} - \mc{T}_r(X)^{\T} \{ R^{\T} \mc{T}_r(X)^{\T} \}^{-} R^{\T} X^{\T}}_F^2 \\
  &= \su \nnorm{(\mc{T}_r(X)^{\T})_{:,i} - \mc{T}_r(X)^{\T} \{ R^{\T} \mc{T}_r(X)^{\T} \}^{-} R^{\T} (X^{\T})_{:,i}}_2^2 \\
  &= \su \nnorm{A (\lambda_i^* - \wt{\lambda}_i)}_2^2 \\
  &= \su \nnorm{\beta_i}_2^2,  \quad \beta_i = A (\lambda_i^* - \wt{\lambda}_i), \; i \in [n].
\end{split} 
\end{align}
Consider the residuals 
\begin{equation}\label{eq:wi}
w_i = b_i - A \lambda_i^* = (X^{\T})_{:,i} - (\mc{T}_r(X)^{\T})_{:,i}.
\end{equation}
By analyzing the structure of (general) sketched regression problems, it can be shown that
\begin{equation}\label{eq:sketchedreg_final_result}
V_r^{\T} R R^{\T} V_r \beta_i = V_r^{\T} R R^{\T} w_i, 
\end{equation}
where $V_r$ is the same matrix as in \eqref{eq:partitioning}. The analysis leading to property \eqref{eq:sketchedreg_final_result} will be given at the end of this proof. In the sequel, we use this property in combination with conditions \textbf{(C1)} and \textbf{(C2)} to deduce the final result. We will first derive a lower bound on the l.h.s.~of \eqref{eq:sketchedreg_final_result} with the help of \textbf{(C2)}, and then we derive an upper bound on the r.h.s.~by
means of \textbf{(C1)}. Combining both, we obtain an upper bound on $\su \nnorm{\beta_i}_2^2$ and in turn on the quantity
$\nnorm{\mc{T}_r(X) - \Pi \mc{T}_r(X)}_F^2$ that we eventually need to bound. 
\vskip1ex
Let $\mc{V}_r \subset \R^d$ denote the column space of $V_r$. Invoking \textbf{(C2)} with $\mc{V} = \mc{V}_r$, the following event
holds with probability at least $1 - \delta_2$:  
\begin{equation*}
\nnorm{R^{\T} V_r v}_2^2 \geq (1 - \eps_2)^2 \nnorm{v}_2^2  \;\, \forall v \in \R^d,
\end{equation*}
or equivalently,
\begin{equation*}
\lambda_{\min}(\Omega) \geq (1- \eps_2)^2, 
\end{equation*}
where $\Omega = V_r^{\T} R R^{\T }V_r$ and $\lambda_{\min}(\cdot)$ denotes the smallest eigenvalue. Conditional on that event, we have that
\begin{align}\label{eq:beta_lowerbound}
\begin{split}
\nnorm{V_r^{\T} R R^{\T} V_r \beta_i}_2^2 &= \beta_i^{\T} \Omega^2 \beta_i \\
                                       &\geq \lambda_{\min}(\Omega^2) \nnorm{\beta_i}_2^2  \\
&\geq (1 - \eps_2)^4 \nnorm{\beta_i}_2^2.
\end{split}
\end{align} 
Next, observe that $V_{:,j}^{\T} w_i = 0$, $j = 1,\ldots,r$, $i=1,\ldots,n$, as follows immediately from the definition of the $\{w_i \}_{i = 1}^n$ in \eqref{eq:wi}. We now apply
\textbf{(C1)} with the following set of vectors:
\begin{equation*}
\mc{S}  = \{ V_{:,j} + \wt{w}_i, \; V_{:,j} - \wt{w}_i, \; i \in [n],\, j \in [r] \},
\end{equation*}
where $\wt{w}_i = w_i / \nnorm{w_i}_2$, $i=1,\ldots,n$. Note that 
$|\mc{S}|  = 2 r n$.   
In the next step, we will establish that with the specified probability, the inner products between 
$V_{:,j}^{\T} w_i$, are preserved up to an additive term of $\eps_1' \nnorm{w_i}_2$, $i \in [n], \; j \in [r]$, where $\eps_1' = \eps_1 / \sqrt{r}$ according to \textbf{(C1)}. 

Recall that for arbitrary $x,y$, it holds that 
$\scp{x}{y} = \frac{1}{4} \left( \nnorm{x + y}_2^2 -  \nnorm{x - y}_2^2  \right)$. 
With $R^{\T}$ being a $(2nr, \eps_1', \delta_1)$ JLT, we therefore have
with probability at least $1 - \delta_1$
\begin{align*}
4 \scp{R^{\T} V_{:,j}}{R^{\T} \wt{w}_i} &= \nnorm{R^{\T} V_{:,j} + R^{\T} \wt{w}_i}_2^2 - \nnorm{R^{\T} V_{:,j} - R^{\T} \wt{w}_i}_2^2 \\
                             &\geq (1 - \eps_1') \nnorm{V_{:,j} + \wt{w}_i}_2^2 - (1 + \eps_1') \nnorm{V_{:,j} - \wt{w}_i}_2^2 \\
                             &= 4 \scp{V_{:,j}}{\wt{w}_i} - 2 \eps_1' \left( \nnorm{V_{:,j}}_2^2 + \nnorm{\wt{w}_i}_2^2 \right) \\
                             &= 4 \scp{V_{:,j}}{\wt{w}_i} - 4 \eps_1'.
\end{align*}
It follows that $\scp{R^{\T} V_{:,j}}{R^{\T} \wt{w}_i} \geq \scp{V_{:,j}}{\wt{w}_i} - \eps_1'$ and in turn also 
$\scp{R^{\T} V_{:,j}}{R^{\T} w_i} \geq \scp{V_{:,j}}{w_i} - \eps_1' \nnorm{w_i}_2$ by homogeneity.

Regarding the upper bound, we argue analogously:
\begin{align*}
4 \scp{R^{\T} V_{:,j}}{R^{\T} \wt{w}_i} &= \nnorm{R^{\T} V_{:,j} + R^{\T} \wt{w}_i}_2^2 - \nnorm{R^{\T} V_{:,j} - R^{\T} \wt{w}_i}_2^2 \\
                             &\leq (1 + \eps_1') \nnorm{V_{:,j} + \wt{w}_i}_2^2 - (1 - \eps_1') \nnorm{V_{:,j} - \wt{w}_i}_2^2 \\
                             &= 4 \scp{V_{:,j}}{\wt{w}_i} + 2 \eps_1' \left( \nnorm{V_{:,j}}_2^2 + \nnorm{\wt{w}_i}_2^2 \right) \\
                             &= 4 \scp{V_{:,j}}{\wt{w}_i} + 4 \eps_1'.
\end{align*}
and thus $\scp{R^{\T} V_{:,j}}{R^{\T} \wt{w}_i} \leq \scp{V_{:,j}}{\wt{w}_i} + \eps_1'$ and in turn $\scp{R^{\T} V_{:,j}}{R^{\T} w_i} \leq \scp{V_{:,j}}{w_i} + \eps_1' \nnorm{w_i}_2$.
We now use these bounds as follows (recall that $\scp{V_{:,j}}{w_i} = 0$, $j \in [r]$, $i \in [n]$): 
\begin{align}\label{eq:wi_lowerbound}
\begin{split}
 \su \nnorm{V_r^{\T} R R^{\T} w_i}_2^2 &= \su \sum_{j = 1}^r \scp{R^{\T} V_{:,j}}{R^{\T} w_i}^2  \\
                                           &\leq \su \sum_{j = 1}^r  (\eps_1')^2\nnorm{w_i}_2^2\\
                                           &= r (\eps_1')^2 \su \nnorm{w_i}_2^2 \\
                                           &= \eps_1^2 \nnorm{X - \mc{T}_r(X)}_F^2
                                         \end{split}
                                       \end{align}
                                       where the last line is immediate from the definition of the $\{ w_i \}_{i = 1}^n$ in \eqref{eq:wi}. Combining \eqref{eq:step2_end}, \eqref{eq:matrixapprox_3},
                                       \eqref{eq:normbetai}, \eqref{eq:sketchedreg_final_result}, \eqref{eq:beta_lowerbound}, \eqref{eq:wi_lowerbound}, we obtain \eqref{eq:mybound_bias_cls} and the assertion of the theorem follows.

                                       \vskip1ex
\noindent In order to finish the proof, it remains to establish \eqref{eq:sketchedreg_final_result} as is done below.   

\vskip2ex 
\noindent For $A \in \R^{d \times n}$, $b \in \R^d$, consider the least squares problem of the form
\begin{equation*}
\min_{\lambda \in \R^n} \nnorm{A\lambda - b}_2^2
\end{equation*}
and the corresponding sketched regression problem with sketching matrix $R^{\T}$
\begin{equation*}
\min_{\lambda} \nnorm{R^{\T} A\lambda - R^{\T} b}_2^2, 
\end{equation*}
Let $\lambda^*$ denote a minimizer of the original least squares problem and let $\wt{\lambda}$ denote the minimizer of the sketched least squares problem. Furthermore,
we write $\mc{U}$ for the matrix of left singular vectors of $A$.  
\vskip1ex
\noindent We then have the following properties:
\begin{itemize} 
\item[(P1)] $A \lambda^* = \mc{U} \alpha$, 
\item[(P2)] $b = A\lambda^* + w$, with $w$ orthogonal to the columns of $\mc{U}$.
\item[(P3)] $A \wt{\lambda} - A \lambda^* = \mc{U} \beta$,
\end{itemize}
for certain vectors $\alpha$ and $\beta$. We now decompose the least squares error when using $\wt{\lambda}$:
\begin{align*}
\nnorm{b - A \wt{\lambda}}_2^2 &=  \nnorm{b - A \lambda^* + A (\lambda^* - \wt{\lambda})}_2^2\\
                         &=  \nnorm{b - A \lambda^*}_2^2 + \nnorm{A (\lambda^* - \wt{\lambda})}_2^2 \\
                        &= \nnorm{w}_2^2 + \nnorm{\mc{U} \beta}_2^2 \\
                        &=  \nnorm{w}_2^2 + \nnorm{\beta}_2^2
\end{align*}
Bringing the sketching matrix $R^{\T}$ into play, we have
\begin{align*}
R^{\T} \mc{U} (\alpha + \beta) &= R^{\T} A \lambda^* + R^{\T} (A\wt{\lambda} - A \lambda^*)\\    
                         &= R^{\T}  A \wt{\lambda} \\
                         &= P_{R^{\T} A} R^{\T} b \\
                         &= P_{R^{\T} \mc{U}} R^{\T} b. 
\end{align*}
Furthermore, we have
\begin{align*}
P_{R^{\T} \mc{U}} R^{\T} b &= P_{R^{\T} \mc{U}} R^{\T} (\mc{U} \alpha + w)\\ 
                   &= R^{\T} \mc{U} \alpha + P_{R^{\T} \mc{U}} R^{\T }w. 
\end{align*}
Combining the previous displays, we obtain that 
\begin{align*}
R^{\T} \mc{U} (\alpha + \beta)  =  R^{\T} \mc{U} \alpha + P_{R^{\T} \mc{U}} R^{\T} w
\end{align*}
and thus 
\begin{equation*}
R^{\T} \mc{U} \beta = P_{R^{\T} \mc{U}} R^{\T} w.
\end{equation*}
Multiplying both sides with $\mc{U}^{\T} R$, this implies 
\begin{align}\label{eq:structure_sketchedleastsquares_final}
\begin{split}
\mc{U}^{\T} R  R^{\T} \mc{U} \beta &=  \mc{U}^{\T} R P_{R^{\T} \mc{U}} R^{\T} w \\
                       &= \mc{U}^{\T} R R^{\T} w. 
                     \end{split}
                   \end{align}
Note that \eqref{eq:structure_sketchedleastsquares_final} has the form as claimed in \eqref{eq:sketchedreg_final_result} with $V_{r}$ playing
the role of $\mc{U}$: according to \eqref{eq:leastsquares_A_b}, this is as it should be since $V_r$ contains the left singular vectors of $\mc{T}_{r}(X)^{\T}$. The proof is thus complete.

\section{Proof of Proposition \ref{prop:randomized_estimation}}
Let us recall that the statement is conditional on $R$, and for what follows only $\{ \omega_l \}_{l = 1}^L$ is considered as random.  
We first verify that $\nnorm{X \omega_l - P_{X_R} X \omega_l}_2^2$ is an unbiased estimator
of $\delta_R^2$, $l \in [L]$. We have 
\begin{align*}
\E[\nnorm{X \omega_l - P_{X_R} X \omega_l}_2^2]  &= \E[\nnorm{(I - P_{X_R}) X \omega_l}_2^2] \\
                                               &= \E[\tr(\omega_l^{\T} X^{\T} (I - P_{X_R}) X \omega_l)] \\
                                               &= \tr(X^{\T} (I - P_{X_R}) X \E[\omega_l \omega_l^{\T}]) \\
                                               &= \tr(X^{\T} (I - P_{X_R}) X) \\
                                               &= \nnorm{X - P_{X_R} X}_F^2. 
\end{align*}
\emph{Concentration.} We now establish concentration for the estimator $\wh{\delta}_R^2$ by invoking results
in \cite{Hsu2012, Laurent2000}. Let $\bm{\omega} \in \R^{d \cdot L}$ be the vector one obtains when stacking 
$\omega_1,\ldots,\omega_L$ vertically. Let us also introduce $\Psi = X^{\T} (I - P_{X_R}) X$ and let $\bm{\Psi} = \frac{1}{L} I_L \otimes \Psi$, where $\otimes$ denotes the Kronecker product. Then $\wh{\delta}_R^2$ can be re-written in the following way:
\begin{align*}
\bm{\omega}^{\T}  \bm{\Psi}   \bm{\omega} &=   \bm{\omega}^{\T} \frac{1}{L}\begin{bmatrix}
\Psi & 0 & \ldots &  \ldots &  0 \\ 
0   & \Psi &  \ldots & \ldots    & 0  \\ 
\vdots  &  0 & \ddots   &  \ddots &   \vdots \\
\vdots &  \vdots  &   \ddots       & \ddots        & 0    \\ 
0      &  0 &  \ldots  & 0   &   \Psi 
\end{bmatrix} \bm{\omega} 
\\
&= \frac{1}{L} \sum_{l = 1}^L \omega_l^{\T} \Psi \omega_l \\
&= \frac{1}{L} \sum_{l = 1}^L \omega_l^{\T} X^{\T} (I - P_{X_R}) X  \omega_l  \\
&= \frac{1}{L} \sum_{l = 1}^L  \nnorm{(I - P_{X_R}) X  \omega_l}_2^2 =\wh{\delta}_R^2.  
\end{align*} 
In other words, $\wh{\delta}_R^2$ can be expressed as a quadratic form in a Gaussian random vector of dimension $dL$ and a positive definite matrix. We can 
thus use the following tail inequalities \cite{Hsu2012, Laurent2000}
\begin{align*}
&\p(\bm{\omega}^{\T} \bm{\Psi} \bm{\omega} > \tr(\bm{\Psi}) + 2 \sqrt{t \tr(\bm{\Psi}^2)} + 2 \nnorm{\bm{\Psi}}_2 t) \leq \exp(-t), \quad t > 0.\\
&\p(\bm{\omega}^{\T} \bm{\Psi} \bm{\omega} < \tr(\bm{\Psi}) - 2 \sqrt{t \tr(\bm{\Psi}^2)}) \leq \exp(-t), \quad t > 0.  %
\end{align*}  
This can be re-written using the following relations:
\begin{align*}
&\tr(\bm{\Psi}) = \tr(\Psi) = \E[\wh{\delta}_R^2] = \delta_R^2, \qquad \sqrt{\tr(\bm{\Psi}^2)} = \nnorm{\bm{\Psi}}_F = \frac{\nnorm{\Psi}_F}{\sqrt{L}} \leq \frac{\tr(\Psi)}{\sqrt{L}},  \\
&\nnorm{\bm{\Psi}}_2 \leq \nnorm{\bm{\Psi}}_F  \leq \tr(\bm{\Psi}),
\end{align*}
\begin{align*}
&\p \left(\wh{\delta}_R^2 >  \delta_R^2 \left(1 + \frac{2 (t + \sqrt{t})}{\sqrt{L}} \right) \right) \leq \exp(-t), \\
&\p \left(\wh{\delta}_R^2 <  \delta_R^2 \left(1 - \frac{2 \sqrt{t}}{\sqrt{L}} \right) \right) \leq \exp(-t), %
\end{align*}
Setting $t = 4$
\begin{align*}
\p \left( \left(1 - \frac{4}{\sqrt{L}} \right) \delta_R^2 \leq \wh{\delta}_R^2 \leq  \delta_R^2 \left(1 + \frac{12}{\sqrt{L}} \right) \right)  \geq 1 - 2 \exp(-4) \geq 0.96.  
\end{align*}
As a result, for any $0 < c < 1$ and any $C > 1$, as long as 
\begin{equation*}
L \geq \max\left\{\frac{16}{(1 - c)^2},  \frac{144}{(C - 1)^2}\right\}
\end{equation*}
it holds that 
\begin{equation*}
\p \left( c \delta_R^2 \leq \wh{\delta}_R^2 \leq  C \delta_R^2 \right)  \geq 1 - 2 \exp(-4) \geq 0.96.
\end{equation*}

\section{Proof of Proposition \ref{prop:columsampling_vs_rp}}
We start with a basic observation to be used several times. Let $\wt{R}$ be an $(d \wedge n) \times k$ random matrix with 
$N(0,1)$ entries. From the rotational invariance of the Gaussian distribution, we have that 
\begin{equation}\label{eq:rotational_invariance}
 V^{\T} R \overset{\mc{D}} = \wt{R}  
\end{equation}
where $\overset{\mc{D}}{=}$ denotes equality in distribution. Turning to property i), using that $\Sigma = \sqrt{n} I$ we have $X_R = \sqrt{n} U V^{\T} R \overset{\mc{D}}{=} \sqrt{n} U \wt{R} = $ according to \eqref{eq:rotational_invariance}. We then compute 
\begin{equation*}
P_{X_R} = X_R (X_R^{\T} X_R)^{-1} X_R^{\T} \overset{\mc{D}}{=} U \wt{R} (\wt{R}^{\T} \wt{R})^{-1} \wt{R}^{\T} U = U P_{\wt{R}} U^{\T},
\end{equation*}
where the inverse exists with probability one. Accordingly,
\begin{align}\label{eq:bias_ortho_exp}
\E \left[\nnorm{X w^* - P_{X_{\wt{R}}} X w^*}_2^2/n \right] 
&= \E[\nnorm{U V^{\T} - U P_{\wt{R}} U^{\T} X}_2^2]  \notag\\
                                         &= \E[\nnorm{U (I - P_{\wt{R}}) V^{\T} w^*}_2^2] \notag\\
                                         &= \E[(\alpha^*)^{\T} (I - P_{\wt{R}}) \alpha^*],
\end{align}
With $\alpha_u^* = \alpha^* / \nnorm{\alpha^*}_2$ and $U_{\wt{R}^{\perp}}$ as a matrix containing a set
of orthonormal basis vectors of $\text{range}(\wt{R})^{\perp}$ as its columns, we have
\begin{align}\label{eq:bias_ortho_exp_2}
\E[(\alpha^*)^{\T} (I - P_{\wt{R}}) \alpha^*] = \nnorm{\alpha^*}_2^2  \E[(\alpha_u^*)^{\T} U_{\wt{R}^{\perp}} U_{\wt{R}^{\perp}}^{\T} (\alpha_u^*)],
\end{align}
Since $\wt{R}$ is Gaussian, $\text{range}(U_{\wt{R}^{\perp}}) \sim \text{Unif}(\textsf{G}(d, d-k))$\footnote{We recall that $\textsf{G}(m, l)$ denotes the set of $l$-dimensional subspaces of $\R^m$.} \cite{James1954}. By rotational invariance, 
\begin{align}\label{eq:bias_ortho_exp_3}
(\alpha_u^*)^{\T} U_{\wt{R}^{\perp}} U_{\wt{R}^{\perp}}^{\T} (\alpha_u^*) = \nnorm{U_{\wt{R}^{\perp}}^{\T} (\alpha_u^*)}_2^2 \overset{\mc{D}}{=}  \nnorm{E_{d - k} u}_2^2 = 1 - k/d,
\end{align} 
where $E_{d-k} \in \R^{(d - k) \times d}$ contains the first $d - k$ canonical basis vectors as its rows and $u \sim \text{Unif}(\mathbb{S}^{d-1})$. Combining
\eqref{eq:bias_ortho_exp}, \eqref{eq:bias_ortho_exp_2} and \eqref{eq:bias_ortho_exp_3} concludes the derivation of the bias. The expression for $\E[\mc{E}(R)]$ given
in the proposition is obtained by adding the variance term $\sigma^2 k/n$. Similarly, we evaluate $\E[\mc{E}(S)]$ by computing its bias. Expanding $P_{X_S}$, we get that
\begin{align}\label{eq:bias_cs_exp_1}
P_{X_S} &= n U V^{\T} S (S^{\T} S)^{-1} S^{\T} V U^{\T} = n U V^{\T} S S^{\T} V U^{\T},        
\end{align}
where we have used that $V V^{\T} = I$ and $S^{\T} S = I$, where the latter property results from the fact that column sampling is done without replacement. It remains to evaluate $\E[S S^{\T}]$. The entries of $S S^{\T}$ are given $(\nscp{S_{i,:}}{S_{j,:}})_{1 \leq i,j \leq d}$, where $S_{l,:}$ denotes the $l$-th row of $S$, $l \in [d]$. We have 
\begin{equation*}
  \E[\nscp{S_{i,:}}{S_{j,:}}] = \E \left[\sum_{l = 1}^k S_{il} S_{jl}  \right] = \sum_{l = 1}^k \p( S_{il} = 1, S_{jl} = 1)
  = \begin{cases}
    \frac{k}{d}& \quad \text{if} \; i = j,               \\
    0  &\quad \text{if} \; i \neq j.
  \end{cases}
\end{equation*}
Putting together the pieces, we obtain that
\begin{equation*}
\E[\nnorm{(I - P_{X_S}) X w^*}_2^2 / n]  = (1 - k/d) \nnorm{\alpha^*}_2^2. 
\end{equation*}
Turning to property ii), the arguments for $\mc{E}(R)$ parallel those used for i), with the difference that $\Sigma = \sqrt{d} I_n$. The subsequent steps are as in i) and are thus omitted. The situation is different
for $\mc{E}(S)$ because the expansion \eqref{eq:bias_cs_exp_1} is no longer valid since $V V^{\T} \neq I$ as $n < d$. Consider a matrix $X$ with dimensions $n < d$ whose matrix of right singular
vectors $V \in \R^{d \times n}$ takes the form
\begin{equation*}
  \left( \begin{array}{l}
    I_n  \\
   \textsf{O}_{n-d, d}
  \end{array} \right),
\end{equation*}
where $\textsf{O}_{n-d, d}$ denotes an $(n-d) \times d$ matrix of zeroes. Note that we still have $V^{\T} V = I_n$, hence this is a valid choice. With this specific form for $V$, we obtain that
\begin{align*}
  P_{X_S} &= U V^{\T} S (S^{\T} V V^{\T} S)^{-1} S^{\T} V U^{\T} \\
            &= U S_{1:n,:} ([S_{1:n,:}]^{\T} S_{1:n,:})^{-1} (S_{1:n,:})^{\T} U^{\T} = U S_{1:n,:} (S_{1:n,:})^{\T} U^{\T},
\end{align*}
where $S_{1:n,:}$ denotes the submatrix of $S$ consisting of its first $n$ rows. It follows that
\begin{align*}
  \E \left[\frac{1}{n} \nnorm{(I - P_{X_S}) Xw^*}_2^2 \right] = \frac{d}{n} (\alpha^*)^{\T} (I - \E[S_{1:n,:} (S_{1:n,:})^{\T} ]) \alpha^* = \frac{d}{n} \nnorm{\alpha^*}_2^2 \left(1 - \frac{k}{d} \right).
\end{align*} 
Regarding property iii), observe that by the rotational invariance according to \eqref{eq:rotational_invariance}
\begin{equation*}
X R = U \Sigma V^{\T} R \overset{\mc{D}}{=} U \Sigma \wt{R}. 
\end{equation*}
Moreover, $X S = U \Sigma (V^{\T} S)$ and $\text{range}(\wt{R}) \overset{\mc{D}}{=} \text{range}(V^{\T} S) \sim \text{Unif}(\textsf{G}(\dwn, k))$ according to \cite{James1954} since the entries of both $R$ and $X$ are i.i.d.~zero-mean Gaussian, thus $\text{range}(XS) \overset{\mc{D}}{=} \text{range}(XR)$.


\section{Proof of Proposition \ref{prop:averaging}}

For property i), observe that the map $A \mapsto \phi(A) \coloneq \nnorm{(I - A) Xw^*}_2^2/n$ from $\R^{n \times n}$ to $\R_+$ is convex, hence
$\phi\left(\frac{1}{B} \sum_{b = 1}^B P_{XR_b} \right) \leq \frac{1}{B}  \sum_{b = 1}^B \phi(P_{XR_b})$. Taking expectations then yields the assertion. Likewise, regarding property
iii), we have
\begin{align*}
  \E\left[\nnorm{\textstyle\frac{1}{B}\textstyle\sum_{b = 1}^B P_{XR_b} \xi}_2^2 \;\Big| \{ P_{XR_b} \}_{b = 1}^B \right]
  &= \E\left[\textstyle\frac{1}{B^2} \textstyle\sum_{b=1}^B \textstyle\sum_{b' = 1}^B \xi^{\T} P_{XR_b} P_{XR_{b'}} \xi  \Big| \{ P_{XR_b} \}_{b = 1}^B  \right]
\end{align*}
The claim then follows by noting that for any pair $(b,b')$, we have $\tr( P_{XR_b} P_{XR_{b'}} ) \leq \nnorm{P_{XR_b}}_F \nnorm{P_{XR_{b'}}}_F = k$.

We finally turn to properties ii) and iv). Consider the operator $\mc{P}_k = \E[P_{XR}]$. We first show that $\text{range}(\mc{P}_k) = \text{range}(X)$. The
inclusion $\text{range}(\mc{P}_k) \subseteq \text{range}(X)$ holds trivially. For the other direction, since $\mc{P}_k$ is symmetric positive definite, it suffices to show that 
$v^{\T} \mc{P}_k  v > 0$ $\forall v \in \text{range}(X)$. Suppose by contradiction that there exists $v \in \text{range}(X)$ s.t. 
\begin{equation*}
v^{\T} \mc{P}_k  v = v^{\T} \E\nolimits[P_{XR}] v = \E\nolimits_R[\nnorm{P_{XR} v}_2^2] = 0,
\end{equation*}
which would imply that $v$ is contained in the orthogonal complement of $\text{range}(XR)$ with probability one, i.e., $v \in \text{null}((XR)^{\T}) \; \Leftrightarrow \, R^{\T} X^{\T} v = 0$ with probability one. This contradicts the fact that the entries of $R$ are from a distribution that is absolutely continuous with respect to the Lebesgue measure. In particular, the fact
that $\text{range}(\mc{P}_k) = \text{range}(X)$ implies that $\mc{P}_k$ has exactly $\dwn$ positive eigenvalues $\{ \eta_j \}_{j = 1}^{\dwn}$ contained in the simplex $\Delta(k) = \{z:\; \textstyle\sum_{j = 1}^{d \wedge n} z_j = k, \; 0 \leq z_j  \leq 1 \}$, noting that $\mc{P}_k$ is an expectation over orthogonal projections onto $k$-dimensional subspaces. The last property of $\mc{P}_k$ to be established
in order to arrive at ii) and iv) is the fact that $U^{\T} \mc{P}_k U = \text{diag}(\eta_1, \ldots, \eta_{\dwn})$, where $U$ is the matrix of left singular vectors of $X$ as its columns. We have
\begin{align*}
\mc{P}_k = \E[P_{XR}] &= \E[X R (R^{\T} X^{\T} X R)^{-1} R^{\T} X^{\T}] \\
          &= \E[U \Sigma V^{\T} R (R^{\T} V \Sigma^2 V^{\T} R)^{-1} R^{\T} V \Sigma U^{\T}] \\
&\hspace*{-18.1ex}\Rightarrow \; U^{\T} \underbrace{\E[P_{XR}]}_{\mc{P}_k} U =  \E[\Sigma V^{\T} \wt{R} (\wt{R}^{\T} V \Sigma^2 V^{\T} \wt{R})^{-1} \wt{R}^{\T} V \Sigma] = \E[\Sigma \wt{R} (\wt{R}^{\T} \Sigma^2 \wt{R})^{-1} \wt{R}^{\T} \Sigma],
\end{align*}
where the last identity uses that $V^{\T} R \overset{\mc{D}}{=} \wt{R}$ by rotational invariance \eqref{eq:rotational_invariance}. It remains to show that the matrix 
\begin{equation}\label{eq:UPkUt}
\E[\Sigma \wt{R} (\wt{R}^{\T} \Sigma^2 \wt{R})^{-1} \wt{R}^{\T} \Sigma]
\end{equation}
is diagonal. This has been shown in \cite{Marzetta2011}, noting that a matrix $A$ is diagonal if and only if $D A D = A$ for all diagonal matrices $D$ with
diagonal elements $\pm 1$; the claim then follows from the fact that 
$D \Sigma \wt{R} \overset{\mc{D}}{=} \Sigma \wt{R}$ and that 
$\wt{R}^{\T} D \Sigma^2 D \wt{R} = \wt{R}^{\T} \Sigma^2 \wt{R}$ for all such $D$. We note that the diagonal elements $\{ \eta_j \}_{j = 1}^{\dwn}$ of the diagonal matrix \eqref{eq:UPkUt} depend only on the singular values $\{ \sigma_j \}_{j = 1}^{\dwn}$ but not on $U$ or $V$ as follows again from rotational invariance. Equipped with the property
$U^{\T} \mc{P}_k U = \text{diag}(\eta_1, \ldots, \eta_{\dwn})$, we compute
\begin{align*}
  \E \left[\nnorm{X w^* - P_{X_R} X w^*}_2^2/n \right] &= \textstyle\frac{1}{n} (w^*)^{\T} X^{\T} \E[I - P_{XR}] X w^* \\
                                                       &= \textstyle\frac{1}{n} (w^*)^{\T} X^{\T} (I - \mc{P}_k) X w^* \\                                                        &= \textstyle\frac{1}{n} \sum_{j = 1}^{d \wedge n} \sigma_j^2 \{\alpha_j^* \}^2 (1 - \eta_j),
\end{align*}
after expanding $X$ in its singular value composition and recalling that $\alpha^* = V^{\T} w^*$. In the same vein, we obtain that
\begin{align*}
  \nnorm{X w^* - \mc{P}_k X w^*}_2^2/n  &= \textstyle\frac{1}{n} \left( \nnorm{Xw^*}_2^2  - 2 (w^*)^{\T} X^{\T} \mc{P}_k X w^* +  (w^*)^{\T} X^{\T} \mc{P}_k^2 X w^* \right)\\
&= \textstyle\frac{1}{n}  \left(  \sum_{j = 1}^{d \wedge n} \sigma_j^2 \{\alpha_j^* \}^2 - 2 
\sum_{j = 1}^{d \wedge n} \sigma_j^2 \{\alpha_j^* \}^2 \eta_j + 
\sum_{j = 1}^{d \wedge n} \sigma_j^2 \{\alpha_j^* \}^2 \eta_j^2 \right) \\
&= \textstyle\frac{1}{n}  \sum_{j = 1}^{d \wedge n} \sigma_j^2 \{\alpha_j^* \}^2 (1 - \eta_j)^2. 
\end{align*}
From $U^{\T} \mc{P}_k^2 U = \text{diag}(\eta_1^2, \ldots, \eta_{\dwn}^2)$, we immediately obtain the first identity in property iv). For the second identity, we let $R'$ be an i.i.d.~copy of $R$ and note that 
\begin{align*}
\tr(\mc{P}_k^2) = \tr(\E[P_{XR}] \E[P_{XR'}]) = \tr(\E[P_{XR} P_{XR'}]) &= \E[\tr(P_{XR} P_{XR'})] \\
&= \E\left[\textstyle\sum_{\ell = 1}^k \cos^2 \theta_{\ell}(\text{range}(XR), \text{range}(XR')) \right],   
\end{align*}
where the last identity is obtained directly from the definition of canonical angles between subspaces \cite{GolubvanLoan}.  


%

\end{document}